\documentclass[11pt]{article}
\usepackage{mathrsfs}
\usepackage{amsmath,amscd,amstext,amsfonts,verbatim}
\usepackage{amsfonts}
\usepackage{bbm}
\usepackage{amscd}
\usepackage{xypic}
\usepackage{amsmath}
\usepackage{amssymb}
\usepackage{amsthm}
\usepackage{bbm}
\usepackage{CJK}
\usepackage{fancyhdr}
\usepackage{hyperref}
\usepackage{indentfirst}
\usepackage{latexsym}
\usepackage{mathrsfs}
 \allowdisplaybreaks 
\usepackage{color}

\usepackage[top=1in,bottom=1in,left=1.25in,right=1.25in]{geometry}
\def\<{\langle}
\def\>{\rangle}
\def\a{\alpha}
\def\b{\beta}

\def\c{\cdot}

\def\D{\Delta}

\def\g{\gamma}

\def\r{\rho}

\def\o{\otimes}

\def\v{\varepsilon}

\def\<{\langle}
\def\>{\rangle}

\date{}
\begin{document}
\renewcommand{\baselinestretch}{1.2}
\renewcommand{\arraystretch}{1.0}
\title{\bf Associative $H$-pseudoalgebras with a semigroup}
\date{}
\author{{\bf Linlin Liu$^{1}$\footnote{Corresponding author. E-mail: liulinlin2016@163.com (Liu), huihuizhengmail@126.com (Zheng)}, Huihui Zheng$^{2}$}\\
{\small  1. School of Science, Henan Institute of Technology, Xinxiang 453003, China}\\
 {\small  2. School of Mathematics and Information Science, Henan Normal University, Xinxiang 453007, China}}

 \maketitle
\begin{center}
\begin{minipage}{14.cm}
\begin{center}{\bf ABSTRACT}\end{center}

Family algebraic structures indexed by a semigroup arise naturally in renormalizations of quantum field theory.
In this paper, we first define the notion of $\Omega$-associative $H$-pseudoalgebra, where the operations are indexed by pairs of elements from a semigroup $\Omega$. Then we construct $\Omega$-associative $H$-pseudoalgebras from associative $H$-pseudoalgebras, $\Omega$-associative algebras, Rota-Baxter family algebras, $\Omega$-type $H$-pseudoalgebras and family-type $H$-pseudoalgebras. Moreover, we investigate the cohomology of $\Omega$-associative $H$-pseudoalgebras and establish that it both induces the cohomology of pseudo-$\mathcal{O}$-operator families and governs the associated formal deformations.
As an application, we show that the first-order deformation of a commutative $\Omega$-associative $H$-pseudoalgebra yields an $\Omega$-Poisson $H$-pseudoalgebra.

 \vskip 0.5cm

{\bf Keywords:}  associative $H$-pseudoalgebra; $\Omega$-associative algebra; cohomology; deformation.
 \vskip 0.5cm
{\bf Mathematics Subject Classification 2020:} 16W99; 16E40; 16S80.
\end{minipage}
\end{center}

\section*{1.  Introduction}
\def\theequation{1. \arabic{equation}}
\setcounter{equation} {0} \hskip\parindent

The notion of conformal algebras was introduced by Kac \cite{K1} as an algebraic language describing the singular part of operator product expansion
in conformal field theory. It has proven useful in the study of vertex algebras, which originated in mathematical physics.
 Roughly speaking, Lie conformal algebras correspond to vertex algebras in the same way that Lie algebras correspond to universal enveloping algebras.

Using the language of pseudotensor categories, a natural generalization of conformal
algebras was introduced in \cite{BD1}, under the name of pseudoalgebras. More precisely, for a cocommutative Hopf algebra $H$ over the
field $k$, we consider the pseudotensor category $\mathcal{M}^{*}(H)$. This category has the same objects as $\mathcal{M}^{l}(H)$ (the category of left $H$-modules) but is equipped with different polylinear maps
$$
Lin_{I}(\{L_i\}, M)=\hbox{Hom}_{H^{\o I}}(\boxtimes_{i\in I}L_i, H^{\o I}\o_H M),
$$
where $I$ is any finite non-empty set and $\boxtimes_{i\in I}: \mathcal{M}^{l}(H)^{\o I}\rightarrow \mathcal{M}^{l}(H^{\o I})$ is the tensor product functor.
An algebra in this pseudotensor category is called an $H$-pseudoalgebra or simply a pseudoalgebra. When $H=\mathbb{C}[\partial]$ (the polynomial algebra), it is actually a conformal algebra, and when $H=k$, it is nothing else but the usual algebra.
In \cite{BD1}, Bakalov, D'Andrea and Kac studied Lie and associative algebras in the pseudotensor category $\mathcal{M}^{*}(H)$, which are called Lie $H$-pseudoalgebras and associative $H$-pseudoalgebras. These structures are closely connected to Hamiltonian formalism in the theory of nonlinear evolution equations \cite{Do}. Later they completed the classification of all irreducible representations of finite simple Lie $H$-pseudoalgebras in a series of papers \cite{BD2,BD3,BD4}.
Furthermore, various other algebraic structures have been introduced in $\mathcal{M}^{*}(H)$, including unital associative pseudoalgebras\cite{Re}, Leibniz $H$-pseudoalgebras \cite{WZ2}, and pre-Lie $H$-pseudoalgebras \cite{LW, WZ1}.

Algebraic structures often appear in `family versions', in which the usual operations are replaced by a family
of operations indexed by a semigroup. The first such example emerged in 2007 in work by Ebrahimi-Fard, Gracia-Bondia, and Patras \cite{Eb} on algebraic aspects of
renormalization in quantum field theory, which Guo later formalized as the Rota-Baxter family algebra \cite{Guo}.
Since then, many scholars begin to pay attention to families algebraic structures \cite{Da, Fo, Ma, ZG, Z, ZM}.
In particular, Aguiar\cite{Ag} introduced $S$-relative algebras (such as $S$-relative associative algebras, also named $\Omega$-associative algebras), a concept that is closely related to family algebras, and he established connections between these types of algebras.
Motivated by this, we study the family version of pseudoalgebra, which leads to the concept of an $\Omega$-associative $H$-pseudoalgebra, a key object of study in this paper.

Deformation problems arise in diverse areas of mathematics, particularly in algebra, algebraic and analytic geometry, and mathematical physics.
The deformation theory was first introduced in \cite{KS} to study complex structures on higher-dimensional manifolds. It was later extended by Gerstenhaber \cite{Ge} to rings and algebras, and by Nijenhuis and Richardson \cite{Ni} to Lie algebras. Their work related deformations of associative algebras to Hochschild cohomology, and those of Lie algebras to Chevalley-Eilenberg cohomology. Recently, cohomology and deformation theory of Rota-Baxter algebras, (twisted) $\mathcal{O}$-operators (families), and $\Omega$-associative conformal algebras have been extensively studied in \cite{Das1, Das2, Da, Tang}. In the present paper, we define the cohomology theory of $\Omega$-associative $H$-pseudoalgebras and associate it with the cohomology of $\mathcal{O}$-operator families of associative $H$-pseudoalgebras (which we call pseudo-$\mathcal{O}$-operator families). We also define the one-parameter formal deformations of $\Omega$-associative $H$-pseudoalgebras and interpret them via lower degree cohomology groups. Furthermore, we show that the first-order deformation of a commutative $\Omega$-associative $H$-pseudoalgebra yields an $\Omega$-Poisson $H$-pseudoalgebra.

This paper is organized as follows. In Section 2, we introduce some basic concepts of associative $H$-pseudoalgebras, dendriform $H$-pseudoalgebras and pseudo-$\mathcal{O}$-operator families. Then we establish the relationships between dendriform $H$-pseudoalgebras and other algebraic structures, namely associative $H$-pseudoalgebras, pre-Lie $H$-pseudoalgebras, as well as pseudo-$\mathcal{O}$-operators.
In Section 3, we investigate two classes of $H$-pseudoalgebras with a semigroup $\Omega$: $\Omega$-type $H$-pseudoalgebras (with operations indexed by pairs of elements from $\Omega$) and family-type $H$-pseudoalgebras (with operations indexed by single elements). Furthermore, we construct $\Omega$-associative $H$-pseudoalgebras from various perspectives, including $\Omega$-associative algebras, associative $H$-pseudoalgebras, Rota-Baxter family algebras, and the aforementioned $\Omega$-type and family-type $H$-pseudoalgebras.
In Section 4, we define the cohomology theory for $\Omega$-associative $H$-pseudoalgebras, which induces a corresponding cohomology for pseudo-$\mathcal{O}$-operator families.
This framework is then applied in Section 5 to study the formal deformations of these pseudoalgebras via lower degree cohomology groups. Finally, in Section 6, we demonstrate that the first-order element of a deformation gives rise to an $\Omega$-Poisson $H$-pseudoalgebra.

Throughout this paper, $k$ is a fixed algebraically closed field of characteristic zero, and $H$ is a cocommutative Hopf algebra over $k$. As usual, we adopt Sweedler's notation \cite{S}. For a coalgebra $C$, we write its comultiplication (with the summation sign omitted) as $\Delta(c) = c_{(1)} \otimes c_{(2)}$ for all $c \in C$.
For any vector space $V$, we define the following permutation maps on tensor powers of $V$:
Let $\sigma: V \otimes V \to V \otimes V$ denote the flip map defined by $\sigma(f \otimes g) = g \otimes f$ for all $f, g \in V$.
For higher tensor powers, we let elements of the symmetric group act by permuting the factors. For example, for all $f, g, h \in V$:
$(12)(f \otimes g \otimes h) = g \otimes f \otimes h$,
$(13)(f \otimes g \otimes h) = h \otimes g \otimes f$.
Similarly, we define the actions of $(23)$, $(123)$, $(132)$, and so on.

\section*{2. $\mathcal{O}$-operator families on associative $H$-pseudoalgebras}
\def\theequation{2. \arabic{equation}}
\setcounter{equation} {0} \hskip\parindent

In this section, we first recall some definitions about associative $H$-pseudoalgebras. Then we define $\mathcal{O}$-operator families on associative $H$-pseudoalgebras, which are called pseudo-$\mathcal{O}$-operator families in the following. Finally, we establish the relationships between dendriform $H$-pseudoalgebras and each of the following: associative $H$-pseudoalgebras, pre-Lie $H$-pseudoalgebras, and pseudo-$\mathcal{O}$-operators.\\

{\bf Definition 2.1}(\cite{BD1}) An $H$-pseudoalgebra is a left $H$-module $A$ equipped with a map $*\in \hbox{Hom}_{H^{\o 2}}(A\o A, H^{\o 2}\o_H A)$, called the \textbf{pseudoproduct}. More precisely, the map
$$
*: A\o A\longrightarrow H^{\o 2}\o_H A,~~~a\o b\mapsto a*b
$$
satisfies
$H^{\o 2}$-linearity: for all $a, b\in A$ and $f, g\in H$,
$$
fa*gb=(f\o g\o_H 1)(a*b).
$$
That is, if $a*b=\sum_{i}(f_i\o g_i)\o_H e_i$, then $fa*gb=\sum_{i}(ff_i\o gg_i)\o_H e_i$.
\\

The pseudoproduct $*: A\o A\longrightarrow H^{\o2}\o_H A$ can be naturally expanded to maps
$$
(H^{\o 2}\o_{H}A)\o A\longrightarrow H^{\o 3}\o_H A,~~A\o (H^{\o 2}\o_{H}A)\longrightarrow H^{\o 3}\o_H A
$$
by defining
$$
(F\o_Hx)*y=\sum_{i}(F\o1)(\D\o id)(f_i\o g_i)\o_H e_i,
$$
$$
x*(F\o_Hy)=\sum_{i}(1\o F)(id\o\D)(f_i\o g_i)\o_H e_i,
$$
for all $F\in H^{\o 2}$, where $x*y=\sum_{i}f_i\o g_i\o_H e_i.$

An $H$-pseudoalgebra $(A,*)$ is said to be \textbf{associative} if the identity $(a* b)* c=a*(b*c)$ holds for all $a, b, c\in A$. It is called \textbf{commutative} if $a*b=((12)\o_H id)(b*a)$ holds for all $a,b\in A$.
\\

{\bf Definition 2.2}(\cite{LW2,WZ3}) Let $(A,*)$ be an associative $H$-pseudoalgebra. An \textbf{$A$-bimodule} is a left $H$-module $M$ equipped with two $H^{\o 2}$-linear maps
$$
*_{l}\in \hbox{Hom}_{H^{\o 2}}(A\o M, H^{\o 2}\o_{H}M),\quad *_{r}\in \hbox{Hom}_{H^{\o 2}}(M\o A, H^{\o 2}\o_{H}M),
$$
satisfying the following identities in $H^{\o 3}\o_H M$ for all $a, b\in A$ and $m\in M$:
\begin{eqnarray*}
&&a*_{l}(b*_{l} m)=(a* b)*_{l} m, \quad m*_{r}(a*b)=(m*_{r}a)*b, \quad (a*_{l}m)*_{r}b=a*_{l}(m*_{r}b).
\end{eqnarray*}
Note that $A$ is an $A$-bimodule with the left and right actions given by the pseudoproduct.
\\

{\bf Definition 2.3} Let $(A,*)$ be an associative $H$-pseudoalgebra and $(M, *_{l}, *_{r})$ an $A$-bimodule. A \textbf{pseudo-$\mathcal{O}$-operator family} (on $M$ over $A$) is a collection
$\{T_{\a}:M\rightarrow A\}_{\a\in \Omega}$ of $H$-linear maps satisfying
\begin{eqnarray}
&&T_{\a}(u)*T_{\b}(v)=T_{\a\b}(T_{\a}(u)*_{l}v+u*_{r}T_{\b}(v)),\label{2.1}
\end{eqnarray}
for all $u,v\in M$ and $\a,\b\in \Omega$.\\

{\bf Remark 2.4} (1) A pseudo-$\mathcal{O}$-operator family is exactly an $\mathcal{O}$-operator family \cite{Da} when $H=k$.

(2) When $M=A$ with the adjoint $A$-bimodule structure and $\Omega=\{0\}$, we recover the notion of Rota-Baxter $H$-operator (of weight 0) \cite{LW}.

(3) A \textbf{pseudo-$\mathcal{O}$-operator} $T: M \to A$ can be viewed as a trivial pseudo-$\mathcal{O}$-operator family by defining $T_\alpha = T$ for all $\alpha \in \Omega$.
\\

{\bf Example 2.5} Let $\{T_{\a}: M\rightarrow A\}_{\a\in \Omega}$ be an $\mathcal{O}$-operator family on the associative algebra $A$
with respect to the bimodule $(M, \cdot_{l},\cdot_{r})$. Then $H\o A$ is an associative $H$-pseudoalgebra with the pseudoproduct
$$
(f\o a)*(g\o b)=f\o g\o_H (1\o ab),~\forall~a, b\in A,f,g\in H.
$$
Note that $H\o M$ is a left $H$-module with the action $h(f\o m)=hf\o m$ for
all $m\in M$ and $f,g\in H$. Define
\begin{eqnarray*}
&&\rho_{l}: (H\o A)\o (H\o M)\rightarrow H^{\o 2}\o_H (H\o M)\\
&&\quad\quad\quad(f\o a)\o(g\o m)\mapsto f\o g\o_H(1\o a\cdot_{l}m),
\end{eqnarray*}
and
\begin{eqnarray*}
&&\rho_{r}: (H\o M)\o (H\o A)\rightarrow H^{\o 2}\o_H (H\o M)\\
&&\quad\quad\quad(f\o m)\o(g\o a)\mapsto f\o g\o_H(1\o m\cdot_{r}a).
\end{eqnarray*}
It is easy to prove that $(H\o M,\rho_{l},\rho_{r})$ is an $(H\o A)$-bimodule. Moreover, $\{id\o T_\a: H\o M\rightarrow H\o A\}_{\a\in \Omega}$
is a pseudo-$\mathcal{O}$-operator family.\\

Let $(A,*)$ be an associative $H$-pseudoalgebra. It is easy to prove that $(A\o k\Omega, \widehat{*})$ is an associative $H$-pseudoalgebra with the left $H$-module action and the pseudoproduct
\begin{eqnarray*}
&&h(a\o\a)=ha\o\a,\quad(a\o\a)\widehat{*}(b\o \b)=a*b\o\a\b,
\end{eqnarray*}
for all $a,b\in A$ and $\a,\b\in \Omega$.
Moreover, if $(M, *_{l}, *_{r})$ is an $A$-bimodule then $(M\o k\Omega,\bullet_{l},\bullet_{r})$ is an $A\o k\Omega$-bimodule with the left and right actions
\begin{eqnarray*}
&&(a\o\a)\bullet_{l}(u\o \b)=a*_{l}u\o\a\b,\quad (u\o \b)\bullet_{r}(a\o\a)=u*_{r}a\o\a\b,
\end{eqnarray*}
for all $a\in A, u\in M$ and $\a,\b\in \Omega$. With this assumption, we have the following result.\\

{\bf Proposition 2.6} $\{T_{\a}:M\rightarrow A\}_{\a\in \Omega}$ is a pseudo-$\mathcal{O}$-operator family if and only if $T: M\o k\Omega\rightarrow A\o k\Omega,u\o\a\mapsto T_\a(u)\o\a$ is a pseudo-$\mathcal{O}$-operator.

{\bf Proof.} For all $u,v\in M$ and $\a,\b\in \Omega$, we have
\begin{eqnarray*}
&&T(u\o\a)\widehat{*}T(v\o\b)\\
&=&(T_\a(u)\o\a)\widehat{*}(T_\b(v)\o\b)\\
&=&T_\a(u)*T_\b(v)\o\a\b\\
&=&T_{\a\b}(T_\a(u)*_lv+u*_rT_\b(v))\o\a\b\\
&=&T(T_\a(u)*_lv\o\a\b+u*_rT_\b(v)\o\a\b)\\
&=&T((T_\a(u)\o\a)\bullet_{l}(v\o\b)+(u\o\a)\bullet_{r}(T_\b(v)\o\b))\\
&=&T(T(u\o\a)\bullet_{l}(v\o\b)+(u\o\a)\bullet_{r}T(v\o\b)),
\end{eqnarray*}
as required.
$\hfill \square$
\\

We now introduce the concept of dendriform $H$-pseudoalgebra, which generalizes the notion of dendriform algebra introduced by Loday \cite{Lo} in the study of algebraic $K$-theory.
\\

 {\bf Definition 2.7} A \textbf{dendriform $H$-pseudoalgebra} $(A, \prec, \succ)$ consists of a left $H$-module $A$ equipped with two operations $\prec,\succ\in \hbox{Hom}_{H^{\o 2}}(A\o A,H^{\o2}\o_H A)$ and such that
\begin{eqnarray*}
&&(x\prec y)\prec z=x\prec(y\prec z+y\succ z),\\
&&(x\succ y)\prec z=x\succ(y\prec z),\\
&&(x\prec y+x\succ y)\succ z=x\succ(y\succ z),
\end{eqnarray*}
for all $x,y,z\in A$.
\\

It is well known that dendriform algebras have close relations with associative algebras, pre-Lie algebras and $\mathcal{O}$-operators. Now we generalize these results to the context of pseudoalgebras. Proposition 2.8 and Proposition 2.9 can be verified by a routine computation. The proof of Proposition 2.10 will be provided later when we consider a more general case in the next subsection (see Proposition 3.21).
\\

 {\bf Proposition 2.8} Let $(A, \prec, \succ)$ be a dendriform $H$-pseudoalgebra.
Define
 $$
 x*y=x\prec y+x\succ y,\forall x,y\in A.
$$
Then $(A,*)$ is an associative $H$-pseudoalgebra.\\

{\bf Proposition 2.9} Let $T: M\rightarrow A$ be a pseudo-$\mathcal{O}$-operator. Then $(M,\prec,\succ)$ is a dendriform $H$-pseudoalgebra, where
\begin{eqnarray*}
&&u\prec v=u*T(v),\quad u\succ v=T(u)*v,
\end{eqnarray*}
for all $u,v\in M$.
\\

Recall that an $H$-pseudoalgebra $(A,*)$ is called a \textbf{pre-Lie $H$-pseudoalgebra} \cite{LW} if it satisfies the additional condition
\begin{eqnarray*}
&&(x* y)*z-x*(y* z)=((12)\o_H id)((y* x)*z-y*(x* z)),\forall x, y, z\in A.
\end{eqnarray*}
With this definition, we have the following result.
\\

{\bf Proposition 2.10} Let $(A, \prec, \succ)$ be a dendriform $H$-pseudoalgebra. Define
 $$
 x*y=x\succ y-((12)\o_H id)y\prec x,\forall x,y\in A.
 $$
Then $(A,*)$ is a pre-Lie $H$-pseudoalgebra.

 \section*{3. $\Omega$-associative $H$-pseudoalgebras and related structures}
\def\theequation{3. \arabic{equation}}

\setcounter{equation} {0} \hskip\parindent

In this section, we introduce the definitions of $\Omega$-associative $H$-pseudoalgebras and related $\Omega$-type $H$-pseudoalgebras, whose operations are indexed by pairs of elements from a semigroup $\Omega$. We also consider family-type $H$-pseudoalgebras, such as dendriform family $H$-pseudoalgebras and pre-Lie family $H$-pseudoalgebras, where the operations are indexed by single elements of $\Omega$. Furthermore, we provide several constructions of $\Omega$-associative $H$-pseudoalgebras.\\

{\bf 3.1 associative and dendriform $H$-pseudoalgebras with a semigroup}\\

{\bf Definition 3.1} An \textbf{$\Omega$-associative $H$-pseudoalgebra} is a left $H$-module $A$ together with a collection of $H^{\o 2}$-linear maps (called the pseudoproduct):
$$
 \{*_{\a,\b}: A\o A\rightarrow H^{\o 2}\o_H A,~~~a\o b\mapsto a*_{\a,\b}b\}_{\a,\b\in\Omega}
$$
satisfying
\begin{eqnarray}
(a*_{\a,\b }b)*_{\a\b,\g}c=a*_{\a,\b\g}(b*_{\b,\g}c), \label{3.1}
\end{eqnarray}
for all $a,b,c\in A$ and $\a,\b,\g\in \Omega$.
\\

{\bf Remark 3.2} (1) When the semigroup $\Omega=\{0\}$, an $\Omega$-associative $H$-pseudoalgebra is precisely an associative $H$-pseudoalgebra.

(2) In particular, for the one-dimensional Hopf algebra $H=k$, an $\Omega$-associative $H$-pseudoalgebra
reduces to an $S$-relative associative algebra \cite{Ag}.  This is equivalent to the structure defined as an $\Omega$-associative algebra in \cite{Da}.
\\

{\bf Definition 3.3}  Let $(A,(*_{\a,\b})_{\a,\b\in \Omega})$ and $(A',(*'_{\a,\b})_{\a,\b\in \Omega})$ be two $\Omega$-associative $H$-pseudo-\\
algebras. A morphism $f: A\rightarrow A'$ of $\Omega$-associative $H$-pseudoalgebras consists of a family of $H$-linear maps $f_{\a}:A\rightarrow A'$, one for each $\a\in \Omega$, satisfying
$$
 f_{\a\b}(x*_{\a,\b}y)=f_{\a}(x)*'_{\a,\b}f_{\b}(y),\quad \forall x,y\in A.
$$

The following result finds the connection between associative $H$-pseudoalgebras and $\Omega$-associative $H$-pseudoalgebras.\\

{\bf Proposition 3.4} $(A\o k\Omega,*)$ is an associative $H$-pseudoalgebra if and only if $(A,(*_{\a,\b})_{\a,\b\in \Omega})$ is an $\Omega$-associative $H$-pseudoalgebra, where
\begin{eqnarray*}
&&(x\o\a)*(y\o\b)=x*_{\a,\b}y\o\a\b,
\end{eqnarray*}
for all $x,y\in A$ and $\a,\b\in \Omega$.

{\bf Proof.} $A\o k\Omega$ is a left $H$-module with the action $h\c (x\o \a)=hx\o \a$, for all $h\in H$, $x\in A$ and $\a\in \Omega$ if and only if $A$ is a left $H$-module.

For all $f,g\in H$, $x, y, z\in A$ and $\a,\b,\g\in \Omega$, we get
\begin{eqnarray*}
&&f(x\o\a)*g(y\o\b)=(fx\o\a)*(gy\o\b)=fx*_{\a,\b}gy\o\a\b,
\end{eqnarray*}
and
\begin{eqnarray*}
&&(f\o g\o_H 1\o 1)((x\o\a)*(y\o\b))\\
&=&(f\o g\o_H 1\o 1)(x*_{\a,\b}y\o\a\b)=(f\o g\o_H 1)(x*_{\a,\b}y)\o\a\b.
\end{eqnarray*}
Similarly, we have
\begin{eqnarray*}
&&((x\o\a)*(y\o\b))*(z\o \g)\\
&=&(x*_{\a,\b}y\o\a\b)*(z\o \g)=(x*_{\a,\b}y)*_{\a\b,\g}z\o \a\b\g,
\end{eqnarray*}
and
\begin{eqnarray*}
&&(x\o\a)*((y\o\b)*(z\o \g))\\
&=&(x\o\a)*((y*_{\b,\g}z)\o\b\g)=x*_{\a,\b\g}(y*_{\b,\g}z)\o \a\b\g.
\end{eqnarray*}
By the above, $(A\o k\Omega,*)$ is an associative $H$-pseudoalgebra if and only if the identities $fx*_{\a,\b}gy=(f\o g\o_H 1)(x*_{\a,\b}y)$
and $(x*_{\a,\b}y)*_{\a\b,\g}z=x*_{\a,\b\g}(y*_{\b,\g}z)$ hold for all $f,g\in H, x,y\in A$ and $\a,\b\in \Omega$, that is, if and only if $(A,(*_{\a,\b})_{\a,\b\in \Omega})$ is an $\Omega$-associative $H$-pseudoalgebra.
$\hfill \square$\\

Now we construct $\Omega$-associative $H$-pseudoalgebra from the ordinary $\Omega$-associative algebras.

{\bf Proposition 3.5} Let $(A, \cdot_{\a,\b})$ be an $\Omega$-associative algebra and $H'$ a subbialgebra of the cocommutative Hopf algebra $H$. Then $H\o H'\o A$ is an $\Omega$-associative $H$-pseudoalgebra, where $H$ acts by left multiplication on the first tensor factor, with the following pseudoproduct
$$
(f\o a\o x)*_{\a,\b}(g\o b\o y)=(f\o ga_1)\o_H(1\o ba_2\o x\cdot_{\a,\b}y),
$$
for all $f\o a\o x, g\o b\o y\in H\o H'\o A$ and $\a,\b\in \Omega$.

{\bf Proof.} It is a routine computation and we omit the details.
$\hfill \square$\\

{\bf Remark 3.6} (1) For a ($G$-graded) associative algebra $A$, the corresponding ($G$-graded) associative $H$-pseudoalgebra structure is established in \cite{BD1, WZ1}. The specific case where $H = H'$ and $A$ is the endomorphism algebra over a finite-dimensional vector space is given in \cite{Re}.

(2) Note that when $H'=k$, $H\o H'\o A\simeq H\o A$ is an $\Omega$-associative $H$-pseudoalgebra with the pseudoproduct
$$
(f\o x)*_{\a,\b}(g\o y)=(f\o g)\o_H x\cdot_{\a,\b}y,
$$
for all $f,g\in H$, $x,y\in A$ and $\a,\b\in \Omega$.
\\

Let $\Omega$ be a semigroup and $\lambda \in k$ be given. Recall that a \textbf{Rota-Baxter family algebra of weight $\lambda$} \cite{Guo, ZG} is an associative algebra $A$ equipped with a collection of linear operators $(P_\a)_{\a \in \Omega}$ on $A$ satisfying
\begin{eqnarray*}
&&P_\alpha(a)P_\beta(b)=P_{\alpha\beta}(P_\alpha(a)b+aP_\beta(b)+\lambda ab), \quad \forall a, b\in A, \alpha, \beta \in \Omega.
\end{eqnarray*}

By combining Proposition 2.10 in \cite{Z} with Proposition 2.3 in \cite{Ag}, such an algebra $(A, (P_\alpha)_{\a \in \Omega})$ induces an $\Omega$-associative algebra structure on $A$ defined by
\begin{eqnarray*}
&&x*_{\a,\b}y=P_\a(x)y+xP_\b(y)+\lambda xy, ~ \forall x,y\in A, \a,\b\in \Omega.
\end{eqnarray*}
Consequently, Proposition 3.5 yields the following corollary.
\\

{\bf Corollary 3.7} Let $(A, (P_\a)_{\a \in \Omega})$ be a Rota-Baxter family algebra of weight $\lambda$ and $H'$ a subbialgebra
 of the cocommutative Hopf algebra $H$. Then $H\o H'\o A$ is an $\Omega$-associative $H$-pseudoalgebra with the pseudoproduct
\begin{eqnarray*}
&&(f\o a\o x)*_{\a,\b}(g\o b\o y)=(f\o ga_1)\o_H(1\o ba_2\o P_{\a}(x)y+xP_{\b}(y)+\lambda xy),
\end{eqnarray*}
for all $f, g\in H, a, b\in H'$ and $x, y\in A$.
\\

{\bf Definition 3.8} An \textbf{$\Omega$-dendriform $H$-pseudoalgebra} is a left $H$-module $D$ equipped with two operations $\prec_{\a,\b},\succ_{\a,\b}\in \hbox{Hom}_{H^{\o 2}}(D\o D,H^{\o2}\o_H D)$ for each pair $(\a,\b)\in \Omega^{2}$, satisfying the following axioms for all $x,y,z\in D$ and $\a,\b,\g\in \Omega$:
\begin{eqnarray}
&&(x\prec_{\a,\b}y)\prec_{\a\b,\g}z=x\prec_{\a,\b\g}(y\prec_{\b,\g}z+y\succ_{\b,\g}z),\label{3.2}\\
&&(x\succ_{\a,\b}y)\prec_{\a\b,\g}z=x\succ_{\a,\b\g}(y\prec_{\b,\g}z),\label{3.3}\\
&&(x\prec_{\a,\b}y+x\succ_{\a,\b}y)\succ_{\a\b,\g}z=x\succ_{\a,\b\g}(y\succ_{\b,\g}z).\label{3.4}
\end{eqnarray}

{\bf Remark 3.9} (1) When $\Omega=\{0\}$ is the trivial semigroup, the structure reduces to a dendriform $H$-pseudoalgebra (see Definition 2.1).

(2) When $H=k$, it recovers the notion of an $S$-relative dendriform algebra \cite{Ag}.
\\

Different from the definition of $\Omega$-dendriform $H$-pseudoalgebra, now we introduce the notion of dendriform family
 $H$-pseudoalgebra, which is a generalization of dendriform family algebra \cite{ZG}. For the latter, the operations are indexed by a single element of $\Omega$ rather than by a pair.\\

{\bf Definition 3.10} A \textbf{dendriform family $H$-pseudoalgebra} is a left $H$-module $D$ equipped with two operations $\prec_{\a},\succ_{\a}\in \hbox{Hom}_{H^{\o 2}}(D\o D,H^{\o2}\o_H D)$ for each $\a\in \Omega$, satisfying the following axioms for all $x,y,z\in D$ and $\a,\b\in \Omega$:
\begin{eqnarray}
&&(x\prec_{\a}y)\prec_{\b}z=x\prec_{\a\b}(y\prec_{\b}z+y\succ_{\a}z),\label{3.5}\\
&&(x\succ_{\a}y)\prec_{\b}z=x\succ_{\a}(y\prec_{\b}z),\label{3.6}\\
&&(x\prec_{\b}y+x\succ_{\a}y)\succ_{\a\b}z=x\succ_{\a}(y\succ_{\b}z).\label{3.7}
\end{eqnarray}

The two types of dendriform $H$-pseudoalgebras with a semigroup mentioned above exhibit the following relationship.\\

{\bf Proposition 3.11} Let $D$ be an $\Omega$-dendriform $H$-pseudoalgebra. Suppose the operations $\prec_{\a,\b}$
are independent of $\a$ (i.e., $\prec_{\a,\b}=\prec_{\b}$) and the operations $\succ_{\a,\b}$ are independent of $\b$ (i.e., $\succ_{\a,\b}=\succ_{\a}$). Then $D$ is a dendriform family
$H$-pseudoalgebra, and vice versa.

{\bf Proof.} Under these assumptions, axioms \eqref{3.2}-\eqref{3.4} specialize to \eqref{3.5}-\eqref{3.7}.
$\hfill \square$\\

In the following, we show that an $\Omega$-dendriform $H$-pseudoalgebra (resp., a dendriform family $H$-pseudoalgebra) induces an $\Omega$-associative $H$-pseudoalgebra.
\\

{\bf Proposition 3.12}  Let $(D,(\succ_{\a,\b})_{\a,\b\in \Omega},(\prec_{\a,\b})_{\a,\b\in \Omega})$ be an $\Omega$-dendriform $H$-pseudoalgebra. Defining
\begin{eqnarray*}
&&x*_{\a,\b}y=x\succ_{\a,\b}y+x\prec_{\a,\b}y,\forall~x,y\in D,\a,\b\in \Omega,
\end{eqnarray*}
turns $(D,(*_{\a,\b})_{\a,\b\in \Omega})$ into an $\Omega$-associative $H$-pseudoalgebra.

{\bf Proof.} Adding equations \eqref{3.2}-\eqref{3.4}, one obtains \eqref{3.1}.
$\hfill \square$\\

Combining Proposition 3.11 with Proposition 3.12, we can get the following conclusion directly.

{\bf Corollary 3.13} Let $(D, \prec_{\a}, \succ_{\a})$ be a dendriform family $H$-pseudoalgebra. Defining
\begin{eqnarray*}
&&x*_{\a,\b}y=x\succ_\a y+x\prec_\b y,\forall x,y\in D,\a,\b\in \Omega,
\end{eqnarray*}
turns $D$ into an $\Omega$-associative $H$-pseudoalgebra.\\

{\bf 3.2 Lie and pre-Lie $H$-pseudoalgebras with a semigroup}\\

From now on, we assume that the semigroup $\Omega$ is commutative unless otherwise specified.\\

{\bf Definition 3.14}  An \textbf{$\Omega$-Lie $H$-pseudoalgebra} consists of a left $H$-module $L$ equipped with a collection of $H^{\o 2}$-linear maps (called the pseudobracket):
$$
 \{[\c*_{\a,\b}\c]: L\o L\rightarrow H^{\o 2}\o_H L,~~~a\o b\mapsto [a*_{\a,\b}b]\}_{\a,\b\in\Omega}
$$
and such that
\begin{eqnarray}
&&[x*_{\a,\b}y]=-(\sigma\o \hbox{id})[y*_{\b,\a}x],\label{3.8}\\
&&[[x*_{\a,\b}y]*_{\a\b,\g}z]=[x*_{\a,\b\g}[y*_{\b,\g}z]]-((12)\o_H id)[y*_{\b,\a\g}[x*_{\a,\g}z]],\label{3.9}
\end{eqnarray}
for all $x,y,z\in L$ and $\a,\b,\g\in \Omega$.\\

{\bf Definition 3.15}  An \textbf{$\Omega$-pre-Lie $H$-pseudoalgebra} consists of a left $H$-module $A$ equipped with a collection of $H^{\o 2}$-linear maps
$$
 \{\circ_{\a,\b}: A\o A\rightarrow H^{\o 2}\o_H A,~~~a\o b\mapsto a\circ_{\a,\b} b\}_{\a,\b\in\Omega}
$$
and such that
\begin{eqnarray*}
&&(x\circ_{\a,\b}y)\circ_{\a\b,\g}z-x\circ_{\a,\b\g}(y\circ_{\b,\g}z)=((12)\o_H id)((y\circ_{\b,\a}x)\circ_{\b\a,\g}z-y\circ_{\b,\a\g}(x\circ_{\a,\g}z)),
\end{eqnarray*}
for all $x,y,z\in A$ and $\a,\b\in \Omega$.
\\

Restricting to the case in which the operations $\circ_{\a,\b}$ are independent of the second index $\b$ in Definition 3.15, we
 obtain the following pre-Lie family $H$-pseudoalgebra, which is a generalization of pre-Lie family algebra \cite{ZM}.\\

{\bf Definition 3.16}  A \textbf{pre-Lie family $H$-pseudoalgebra} consists of a left $H$-module $A$ equipped with an $H^{\o 2}$-linear map $\circ_{\a}: A\o A\rightarrow H^{\o 2}\o_H A$ for each $\a\in \Omega$ and such that
\begin{eqnarray*}
&&(x\circ_{\a}y)\circ_{\a\b}z-x\circ_{\a}(y\circ_{\b}z)=((12)\o_H id)((y\circ_{\b}x)\circ_{\b\a}z-y\circ_{\b}(x\circ_{\a}z)),
\end{eqnarray*}
for all $x,y,z\in A$ and $\a,\b\in \Omega$.
\\

{\bf Remark 3.17} (1) An $\Omega$-pre-Lie $H$-pseudoalgebra (resp. a pre-Lie family $H$-pseudoalgebra) becomes a pre-Lie $H$-pseudoalgebra if $\Omega={0}$, and an $\Omega$-pre-Lie algebra \cite{Ag} (resp. pre-Lie family algebra) if $H=k$.

(2) Clearly, any $\Omega$-associative $H$-pseudoalgebra is an $\Omega$-pre-Lie $H$-pseudoalgebra.
\\

{\bf Example 3.18} Let $(A, \cdot_{\a,\b})$ be an $\Omega$-pre-Lie algebra and $C$ a subcoalgebra of $H$. Then $H\o C\o A$ is a left $H$-module with $H$ acting by left multiplication on the first tensor factor.
Moreover, $H\o C\o A$ becomes an $\Omega$-pre-Lie $H$-pseudoalgebra with the pseudoproduct
\begin{eqnarray}
(f\o a\o x)*_{\a,\b}(g\o b\o y)=faS(b_1)\o g\o_H 1\o b_2\o x\cdot_{\a,\b}y,\label{3.10}
\end{eqnarray}
for all $f\o a\o x, g\o b\o y\in H\o C\o A$ and $\a,\b\in \Omega$.

If $A$ is an $\Omega$-associative algebra, then $H\o C\o A$ is an $\Omega$-associative $H$-pseudoalgebra with the above pseudoproduct.
If $A$ is a pre-Lie family algebra, then $H\o C\o A$ is a pre-Lie family $H$-pseudoalgebra with pseudoproduct \eqref{3.10} by letting $\cdot_{\a,\b}=\cdot_{\a}$ and $*_{\a,\b}=*_{\a}$.\\

In \cite{WZ1} the author showed that a $G$-graded pre-Lie $H$-pseudoalgebra induces a $G$-graded Lie $H$-pseudoalgebra. Similarly, we have the following proposition.\\

{\bf Proposition 3.19}  Let $(A, (\circ_{\a,\b})_{\a,\b\in\Omega})$ be an $\Omega$-pre-Lie $H$-pseudoalgebra. Define
\begin{eqnarray}
&&[x*_{\a,\b}y]=x\circ_{\a,\b} y-(\sigma\o_H id)(y\circ_{\b,\a} x),\label{3.11}
\end{eqnarray}
for all $x,y\in A$ and $\a,\b\in \Omega$. Then $(A, ([\cdot*_{\a,\b}\cdot])_{\a,\b\in \Omega})$ is an $\Omega$-Lie $H$-pseudoalgebra.

{\bf Proof.} Condition \eqref{3.8} follows immediately from \eqref{3.11}. For all $x,y,z\in A$ and $\a,\b,\g\in \Omega$, we have
\begin{eqnarray*}
[[x*_{\a,\b}y]*_{\a\b,\g}z]&=&(x\circ_{\a,\b}y)\circ_{\a\b,\g}z-((12)\o_H id)(y\circ_{\b,\a}x)\circ_{\a\b,\g}z\\
&&-((132)\o_H id)z\circ_{\g,\a\b}(x\circ_{\a,\b}y)+((13)\o_H id)z\circ_{\g,\a\b}(y\circ_{\b,\a}x),
\end{eqnarray*}
and
\begin{eqnarray*}
[x*_{\a,\b\g}[y*_{\b,\g}z]]&=&x\circ_{\a,\b\g}(y\circ_{\b,\g}z)-((23)\o_H id)x\circ_{\a,\b\g}(z\circ_{\g,\b}y)\\
&&-((123)\o_H id)(y\circ_{\b,\g}z)\circ_{\b\g,\a}x+((13)\o_H id)(z\circ_{\g,\b}y)\circ_{\b\g,\a}x.
\end{eqnarray*}
Similarly, we have
\begin{eqnarray*}
&&-((12)\o_H id)[y*_{\b,\a\g}[x*_{\a,\g}z]]\\
&=&-((12)\o_H id)y\circ_{\b,\a\g}(x\circ_{\a,\g}z)+((123)\o_H id)y\circ_{\b,\a\g}(z\circ_{\g,\a}x)\\
&&+((23)\o_H id)(x\circ_{\a,\g}z)\circ_{\a\g,\b}y+((132)\o_H id)(z\circ_{\g,\a}x)\circ_{\a\g,\b}y.
\end{eqnarray*}
Now we compute
\begin{eqnarray*}
&&[x*_{\a,\b\g}[y*_{\b,\g}z]]-((12)\o_H id)[y*_{\b,\a\g}[x*_{\a,\g}z]]-[[x*_{\a,\b}y]*_{\a\b,\g}z]\\
&=&x\circ_{\a,\b\g}(y\circ_{\b,\g}z)-(x\circ_{\a,\b}y)\circ_{\a\b,\g}z-((12)\o_H id)y\circ_{\b,\a\g}(x\circ_{\a,\g}z)\\
&&+((12)\o_H id)(y\circ_{\b,\a}x)\circ_{\a\b,\g}z-((23)\o_H id)(x\circ_{\a,\b\g}(z\circ_{\g,\b}y)-(x\circ_{\a,\g}z)\circ_{\a\g,\b}y\\
&&+((12)\o_H id)(z\circ_{\g,\a}x)\circ_{\a\g,\b}y-((12)\o_H id)z\circ_{\g,\a\b}(x\circ_{\a,\b}y))\\
&&-((123)\o_H id)((y\circ_{\b,\g}z)\circ_{\b\g,\a}x-y\circ_{\b,\a\g}(z\circ_{\g,\a}x)-((12)\o_H id)(z\circ_{\g,\b}y)\circ_{\b\g,\a}x\\
&&+((12)\o_H id)z\circ_{\g,\a\b}(y\circ_{\b,\a}x))\\
&=&0.
\end{eqnarray*}
So \eqref{3.9} holds and we complete the proof.\hfill $\square$\\

The following corollaries can be obtained directly.\\

{\bf Corollary 3.20} Let $(A, (*_{\a,\b})_{\a,\b\in\Omega})$ be an $\Omega$-associative $H$-pseudoalgebra. Define
\begin{eqnarray}
&&[x*_{\a,\b}y]=x*_{\a,\b} y-(\sigma\o_H id)(y*_{\b,\a} x),\label{3.12}
\end{eqnarray}
for all $x,y\in A$ and $\a,\b\in \Omega$. Then $(A, ([\cdot*_{\a,\b}\cdot])_{\a,\b\in \Omega})$ is an $\Omega$-Lie $H$-pseudoalgebra.
\\

{\bf Corollary 3.21} Let $(A, \circ_{\a})$ be a pre-Lie family $H$-pseudoalgebra. Define
\begin{eqnarray*}
&&[x*_{\a,\b}y]=x\circ_{\a} y-(\sigma\o id)y\circ_{\b} x,
\end{eqnarray*}
for all $x,y\in A$ and $\a,\b\in \Omega$. Then $(A, ([\cdot*_{\a,\b}\cdot])_{\a,\b\in \Omega})$ is an $\Omega$-Lie $H$-pseudoalgebra.
\\

Now we establish the relationship between $\Omega$-pre-Lie $H$-pseudoalgebras and $\Omega$-dendriform $H$-pseudoalgebras.\\

{\bf Proposition 3.22}  Let $(D,(\succ_{\a,\b})_{\a,\b\in \Omega},(\prec_{\a,\b})_{\a,\b\in \Omega})$ be an $\Omega$-dendriform $H$-pseudoalgebra. Define
\begin{eqnarray*}
&&x\circ_{\a,\b}y=x\succ_{\a,\b}y-(\sigma\o_H id)y\prec_{\b,\a}x,
\end{eqnarray*}
for all $x,y\in D$ and $\a,\b\in \Omega$. Then $(D,(\circ_{\a,\b})_{\a,\b\in \Omega})$ is an $\Omega$-pre-Lie $H$-pseudoalgebra.

{\bf Proof.} For all $x,y\in D$ and $\a,\b\in \Omega$, we have
\begin{eqnarray*}
&&x\circ_{\a,\b\g}(y\circ_{\b,\g}z)\\
&=&x\circ_{\a,\b\g}(y\succ_{\b,\g}z-(\sigma\o_H id)z\prec_{\g,\b}y)\\
&=&x\succ_{\a,\b\g}(y\succ_{\b,\g}z)-((23)\o_H id)x\succ_{\a,\b\g}(z\prec_{\g,\b}y)\\
&&-((123)\o_H id)(y\succ_{\b,\g}z)\prec_{\b\g,\a}x+((13)\o_H id)(z\prec_{\g,\b}y)\prec_{\b\g,\a}x,
\end{eqnarray*}
and
\begin{eqnarray*}
&&(x\circ_{\a,\b}y)\circ_{\a\b,\g}z\\
&=&(x\succ_{\a,\b}y-(\sigma\o_H id)y\prec_{\b,\a}x)\circ_{\a\b,\g}z\\
&=&(x\succ_{\a,\b}y)\succ_{\a\b,\g}z-((12)\o_H id)(y\prec_{\b,\a}x)\succ_{\a\b,\g}z\\
&&-((132)\o_H id)z\prec_{\g,\a\b}(x\succ_{\a,\b}y)+((13)\o_H id)z\prec_{\g,\a\b}(y\prec_{\b,\a}x).\\
\end{eqnarray*}
Now we compute
\begin{eqnarray*}
&&x\circ_{\a,\b\g}(y\circ_{\b,\g}z)-(x\circ_{\a,\b}y)\circ_{\a\b,\g}z\\
&=&x\succ_{\a,\b\g}(y\succ_{\b,\g}z)-(x\succ_{\a,\b}y)\succ_{\a\b,\g}z+((12)\o_H id)(y\prec_{\b,\a}x)\succ_{\a\b,\g}z\\
&&-((23)\o_H id)x\succ_{\a,\b\g}(z\prec_{\g,\b}y)-((123)\o_H id)(y\succ_{\b,\g}z)\prec_{\b\g,\a}x\\
&&+((13)\o_H id)(z\prec_{\g,\b}y)\prec_{\b\g,\a}x+((132)\o_H id)z\prec_{\g,\a\b}(x\succ_{\a,\b}y)\\
&&-((13)\o_H id)z\prec_{\g,\a\b}(y\prec_{\b,\a}x)\\
&=&(x\prec_{\a,\b}y)\succ_{\a\b,\g}z+((12)\o_H id)(y\succ_{\b,\a\g}(x\succ_{\a,\g}z)-(y\succ_{\b,\a}x)\succ_{\a\b,\g}z)\\
&&-((23)\o_H id)(x\succ_{\a,\g}z)\prec_{\a\g,\b}y-((123)\o_H id)y\succ_{\b,\a\g}(z\prec_{\g,\a}x)\\
&&+((13)\o_H id)z\prec_{\g,\a\b}(y\succ_{\b,\a}x)+((132)\o_H id)(z\prec_{\g,\a}x)\prec_{\a\g,\b}y\\
&&-((132)\o_H id)z\prec_{\g,\a\b}(x\prec_{\a,\b}y)\\
&=&((12)(12)\o_H id)(x\prec_{\a,\b}y)\succ_{\a\b,\g}z+((12)\o_H id)y\succ_{\b,\a\g}(x\succ_{\a,\g}z)\\
&&-((12)\o_H id)(y\succ_{\b,\a}x)\succ_{\a\b,\g}z-((12)(123)\o_H id)(x\succ_{\a,\g}z)\prec_{\a\g,\b}y\\
&&-((12)(23)\o_H id)y\succ_{\b,\a\g}(z\prec_{\g,\a}x)+((12)(132)\o_H id)z\prec_{\g,\a\b}(y\succ_{\b,\a}x)\\
&&+((12)(13)\o_H id)((z\prec_{\g,\a}x)\prec_{\a\g,\b}y-z\prec_{\g,\a\b}(x\prec_{\a,\b}y))\\
&=&((12)\o_H id)(y\succ_{\b,\a\g}(x\succ_{\a,\g}z)-((23)\o_H id)y\succ_{\b,\a\g}(z\prec_{\g,\a}x)\\
&&-((123)\o_H id)(x\succ_{\a,\g}z)\prec_{\a\g,\b}y+((13)\o_H id)(z\prec_{\g,\a}x)\prec_{\a\g,\b}y)\\
&&-((12)\o_H id)((y\succ_{\b,\a}x)\succ_{\a\b,\g}z-((12)\o_H id)(x\prec_{\a,\b}y)\succ_{\a\b,\g}z\\
&&-((132)\o_H id)z\prec_{\g,\a\b}(y\succ_{\b,\a}x)+((13)\o_H id)z\prec_{\g,\a\b}(x\prec_{\a,\b}y))\\
&=&((12)\o_H id)(y\circ_{\b,\a\g}(x\circ_{\a,\g}z)-(y\circ_{\b,\a}x)\circ_{\b\a,\g}z),
\end{eqnarray*}
which completes the proof.\hfill $\square$\\

Similarly, for the family case we have the following result.\\

{\bf Proposition 3.23}  Let $(D,(\succ_{\a})_{\a\in \Omega},(\prec_{\a})_{\a\in \Omega})$ be a dendriform family $H$-pseudoalgebra. Define
\begin{eqnarray*}
&&x\circ_{\a}y=x\succ_{\a}y-(\sigma\o_H id)y\prec_{\a}x,
\end{eqnarray*}
for all $x,y\in D$ and $\a\in \Omega$. Then $(D,(\circ_{\a})_{\a\in \Omega})$ is a pre-Lie family $H$-pseudoalgebra.
\\

{\bf 3.3 Zinbiel $H$-pseudoalgebras with a semigroup}\\

{\bf Definition 3.24} An \textbf{$\Omega$-zinbiel $H$-pseudoalgebra} consists of a left $H$-module $A$ equipped with an operation $*_{\a,\b}\in \hbox{Hom}_{H^{\o 2}}(A\o A,H^{\o2}\o_H A)$ for each pair $(\a,\b)\in \Omega^{2}$ and such that
\begin{eqnarray}
&&x*_{\a,\b\g}(y*_{\b,\g}z)=(x*_{\a,\b}y)*_{\a\b,\g}z+((12)\o_H id)(y*_{\b,\a}x)*_{\a\b,\g}z,\label{3.13}
\end{eqnarray}
for all $x,y,z\in A$ and $\a,\b,\g\in \Omega$.\\

 Consider now a left $H$-module $A$ equipped with a map $*_{\a}\in \hbox{Hom}_{H^{\o 2}}(A\o A,H^{\o2}\o_H A)$ for each $\a\in \Omega$ and such that
 \begin{eqnarray*}
&&x*_{\a}(y*_{\b}z)=(x*_{\a}y)*_{\a\b}z+((12)\o_H id)(y*_{\b}x)*_{\a\b}z,~\forall~x,y,z\in A, \a,\b\in \Omega.
\end{eqnarray*}
Then $(A, (*_{\a})_{\a\in \Omega})$ is called a \textbf{zinbiel family $H$-pseudoalgebra}.

Note that a zinbiel family $H$-pseudoalgebra is the same as an $\Omega$-zinbiel $H$-pseudoalgebra in which the operations $*_{\a,\b}$ are independent of $\b$.\\

Now we establish the relationship between $\Omega$-dendriform $H$-pseudoalgebras and $\Omega$-zinbiel $H$-pseudoalgebras. Prior to this, we have the following lemma.\\

{\bf Lemma 3.25} Let $(A, (*_{\a,\b})_{\a,\b\in \Omega})$ be an $\Omega$-zinbiel $H$-pseudoalgebra. Then
\begin{eqnarray}
&&x*_{\a,\b\g}(y*_{\b,\g}z)=((12)\o_H id)y*_{\b,\a\g}(x*_{\a,\g}z),\label{3.14}
\end{eqnarray}
for all $x,y,z\in A$ and $\a,\b,\g\in \Omega$.

{\bf Proof.} Replacing $(x,y)$ for $(y,x)$ and $(\a,\b)$ for $(\b,\a)$ in \eqref{3.13}, we have
\begin{eqnarray*}
&&y*_{\b,\a\g}(x*_{\a,\g}z)\\
&=&(y*_{\b,\a}x)*_{\b\a,\g}z+((12)\o_H id)(x*_{\a,\b}y)*_{\a\b,\g}z\\
&=&((12)\o_H id)((x*_{\a,\b}y)*_{\a\b,\g}z+((12)\o_H id)(y*_{\b,\a}x)*_{\b\a,\g}z)\\
&=&((12)\o_H id)x*_{\a,\b\g}(y*_{\b,\g}z).
\end{eqnarray*}
This completes the proof.\hfill $\square$\\

{\bf Proposition 3.26}  Let $(D,(\succ_{\a,\b})_{\a,\b\in \Omega},(\prec_{\a,\b})_{\a,\b\in \Omega})$ be an $\Omega$-dendriform $H$-pseudoalgebra
and $(A, (*_{\a,\b})_{\a,\b\in \Omega})$ an $\Omega$-zinbiel $H$-pseudoalgebra.

(1) Suppose $x\succ_{\a,\b}y=(\sigma\o id)y\prec_{\b,\a}x$. Then defining $x*_{\a,\b}y=x\succ_{\a,\b}y,\forall x,y\in D,\a,\b\in \Omega$ turns $D$ into an $\Omega$-zinbiel $H$-pseudoalgebra.

(2) Define $x\prec_{\a,\b}y=(\sigma\o_H id)y*_{\b,\a}x$ and $x\succ_{\a,\b}y=x*_{\a,\b}y$. Then $A$ is an $\Omega$-dendriform $H$-pseudoalgebra. Moreover, $x\succ_{\a,\b}y=(\sigma\o_H id)y\prec_{\b,\a}x$.

{\bf Proof.}(1) For all $x, y, z\in D$ and $\a,\b, \g\in \Omega$, we have
\begin{eqnarray*}
&&(x*_{\a,\b}y)*_{\a\b,\g}z+((12)\o_H id)(y*_{\b,\a}x)*_{\b\a,\g}z\\
&=&(x\succ_{\a,\b}y)\succ_{\a\b,\g}z+((12)\o_H id)(y\succ_{\b,\a}x)\circ_{\b\a,\g}z\\
&=&(x\succ_{\a,\b}y)\succ_{\a\b,\g}z+(x\prec_{\a,\b}y)\circ_{\a\b,\g}z\\
&=&x\succ_{\a,\b\g}(y\succ_{\b,\g}z)=x*_{\a,\b\g}(y*_{\b,\g}z),
\end{eqnarray*}
as desired.

(2) For all $x, y, z\in A$ and $\a,\b, \g\in \Omega$, we obtain
\begin{eqnarray*}
&&(x\prec_{\a,\b\g}(y\prec_{\b,\g}z+y\succ_{\b,\g}z)\\
&=&((13)\o_H id)(z*_{\g,\b}y)*_{\b\g,\a}x+((123)\o_H id)(y*_{\b,\g}z)*_{\b\g,\a}x\\
&=&((13)\o_H id)((z*_{\g,\b}y)*_{\b\g,\a}x+((12)\o_H id)(y*_{\b,\g}z)*_{\b\g,\a}x)\\
&=&((13)\o_H id)z*_{\g,\a\b}(y*_{\b,\a}x)\quad(\hbox{by Eq.}~\eqref{3.13})\\
&=&((\sigma\o_H id)y*_{\b,\a}x)\prec_{\a\b,\g}z\\
&=&(x\prec_{\a,\b}y)\prec_{\a\b,\g}z,
\end{eqnarray*}
\begin{eqnarray*}
&&(x\succ_{\a,\b}y)\prec_{\a\b,\g}z\\
&=&((132)\o_H id)z*_{\g,\a\b}(x*_{\a,\b}y)\\
&=&((23)(12)\o_H id)z*_{\g,\a\b}(x*_{\a,\b}y)\\
&=&((23)\o_H id)x*_{\a,\b\g}(z*_{\g,\b}y)\quad(\hbox{by Eq.}~\eqref{3.14})\\
&=&x\succ_{\a,\b\g}((\sigma\o_H id)z*_{\g,\b}y)\\
&=&x\succ_{\a,\b\g}(y\prec_{\b,\g}z),
\end{eqnarray*}
and
\begin{eqnarray*}
&&(x\prec_{\a,\b}y+x\succ_{\a,\b}y)\succ_{\a\b,\g}z\\
&=&((\sigma\o_H id)y*_{\b,\a}x+x*_{\a,\b}y)*_{\a\b,\g}z\\
&=&((12)\o_H id)(y*_{\b,\a}x)*_{\a\b,\g}z+(x*_{\a,\b}y)*_{\a\b,\g}z\\
&=&x*_{\a,\b\g}(y*_{\b,\g}z)\quad(\hbox{by Eq.}~\eqref{3.13})\\
&=&x\succ_{\a,\b\g}(y\succ_{\b,\g}z).
\end{eqnarray*}
So we complete the proof.
\hfill $\square$\\

Restricting to the case in which the operations $\prec_{\a,\b}$ are independent of $\a$, $\succ_{\a,\b}$ are independent of $\b$ and $*_{\a,\b}$ are independent of the second index $\b$.
we obtain corresponding results for the family case.\\

{\bf Proposition 3.27}  Let $(D,(\succ_{\a})_{\a\in \Omega},(\prec_{\a})_{\a\in \Omega})$ be a dendriform family $H$-pseudoalgebra
and $(A, (*_{\a,\b})_{\a,\b\in \Omega})$ a zinbiel family $H$-pseudoalgebra.

(1) Suppose $x\succ_{\a}y=(\sigma\o id)y\prec_{\a}x$. Then defining $x*_{\a}y=x\succ_{\a}y,\forall x,y\in D,\a\in \Omega$ turns $D$ into a zinbiel family $H$-pseudoalgebra.

(2) Define $x\prec_{\a}y=(\sigma\o_H id)y*_{\a}x$ and $x\succ_{\a}y=x*_{\a}y$. Then $A$ is a dendriform family $H$-pseudoalgebra. Moreover, $x\succ_{\a}y=(\sigma\o_H id)y\prec_{\a}x$.
\\

An $\Omega$-associative $H$-pseudoalgebra $(A,(*_{\a,\b})_{\a,\b\in \Omega})$ is said to be \textbf{commutative} if $x*_{\a,\b} y=(\sigma\o_H id)(y*_{\b,\a} x)$ for all $x, y\in A$ and $\a,\b\in \Omega$. Now we construct a commutative $\Omega$-associative $H$-pseudoalgebra from an $\Omega$-zinbiel $H$-pseudoalgebra.\\

{\bf Proposition 3.28}  Let $(A, (*_{\a,\b})_{\a,\b\in \Omega})$ be an $\Omega$-zinbiel $H$-pseudoalgebra. Define
\begin{eqnarray}
&&x\widehat{*}_{\a,\b}y=x*_{\a,\b}y+(\sigma\o_H id)y*_{\b,\a}x,\label{3.15}
\end{eqnarray}
for all $x,y\in A$ and $\a,\b\in \Omega$. Then $(A, (\widehat{*}_{\a,\b})_{\a,\b\in \Omega})$ is a commutative $\Omega$-associative $H$-pseudoalgebra.

{\bf Proof.} By Proposition 3.12 and Proposition 3.26(2), $(A, (\widehat{*}_{\a,\b})_{\a,\b\in \Omega})$ is an $\Omega$-associative $H$-pseudoalgebra. The commutativity follows immediately from \eqref{3.15}, which completes the proof.
\hfill $\square$

\section*{4. Cohomology theory of $\Omega$-associative $H$-pseudoalgebras}
\def\theequation{4. \arabic{equation}}
\setcounter{equation} {0} \hskip\parindent

In this section, we develop the cohomology theory for $\Omega$-associative $H$-pseudoalgebras and extend it to define the cohomology of pseudo-$\mathcal{O}$-operator families.

Throughout this section, we assume that $\Omega$ is a semigroup with unit $1\in \Omega$. The unital condition of $\Omega$ is only used in defining the coboundary
operator at degree 0.\\

{\bf 4.1 The cohomology of $\Omega$-associative $H$-pseudoalgebras}\\

In this subsection, we first define bimodules over $\Omega$-associative $H$-pseudoalgebras and then introduce their cohomology theory.\\

{\bf Definition 4.1}  Let $(A,(*_{\a,\b})_{\a,\b\in \Omega})$ be an $\Omega$-associative $H$-pseudoalgebra. An $A$-bimodule is a left $H$-module $M$ equipped with a collection
\begin{align*}
\begin{Bmatrix}
\r^{l}_{\a,\b}: A\o M\rightarrow H^{\o 2}\o_H M,~~(a, m)\mapsto a\cdot_{\a,\b}m\\
\r^{r}_{\a,\b}: M\o A\rightarrow H^{\o 2}\o_H M,~~(m, a)\mapsto m\cdot_{\a,\b}a
\end{Bmatrix}_{\a,\b\in \Omega}
\end{align*}
of $H^{\o 2}$-linear maps that satisfy the following axioms:
\begin{eqnarray*}
&&(a*_{\a,\b}b)\cdot_{\a\b,\g}m=a\c_{\a,\b\g}(b\c_{\b,\g}m),\\
&&(a\c_{\a,\b}m)\cdot_{\a\b,\g}b=a\c_{\a,\b\g}(m\c_{\b,\g}b),\\
&&(m\c_{\a,\b}a)\cdot_{\a\b,\g}b=m\c_{\a,\b\g}(a*_{\b,\g}b),
\end{eqnarray*}
for all $a,b\in A$, $m\in M$ and $\a, \b, \g\in \Omega$.

Let $(A,(*_{\a,\b})_{\a,\b\in \Omega})$ be an $\Omega$-associative $H$-pseudoalgebra and $M$ an $A$-bimodule. Now we consider the Hochschild cohomology complex $C^{\bullet}(A,M)$ for an $\Omega$-associative $H$-pseudoalgebra $A$ with coefficients in an $A$-bimodule.

For each $n> 0$, the space of $n$-cochains $C^{n}_{\Omega}(A, M)$ consists of all maps
$$
f=\{f_{\a_1,\cdots, \a_n}\}_{\a_1,\cdots, \a_n\in \Omega}\in \hbox{Hom}_{H^{\o n}}(A^{\o n},H^{\o n}\o_H M).
$$
More precisely, $f$ satisfies
$$
f_{\a_1,\cdots, \a_n}(h_1a_1,\cdots, h_na_n)=(h_1\o\cdots\o h_n\o 1)f_{\a_1,\cdots, \a_n}(a_1,\cdots, a_n).
$$
For $n=0$, $C^{0}_{\Omega}(A, M)=k\o_H M\simeq M/H_+M$, where $H_+=\{h\in H|\varepsilon(h)=0\}$ is the augmentation ideal of $H$.

The differential $d_{\Omega}^{n}: C^{n}_{\Omega}(A, M)\rightarrow C^{n+1}_{\Omega}(A, M)$ is given by

\begin{eqnarray*}
&&d_{\Omega}^{0}(1\o_Hm)_{\a}(a)=\sum_i(id\o\v)(h_{i})m_i-\sum_j(\v\o id)(f_{j})n_j,
\end{eqnarray*}
where $a\c_{\a,1}m=\sum_i h_{i}\o_H m_{i}$ and $m\c_{1,\a}a=\sum_j f_{j}\o_H n_{j}$.

\begin{eqnarray}
\nonumber&&(d_{\Omega}^{n\geq1}f)_{\a_1,\cdots,\a_{n+1}}(a_1,\cdots, a_{n+1})\\
\nonumber&=&a_1\cdot_{\a_1,\a_2\cdots \a_{n+1}}f_{\a_2,\cdots, \a_{n+1}}(a_2,\cdots, a_{n+1})\\
\nonumber&&+\sum_{i=1}^{n}(-1)^{i}f_{\a_1,\cdots, \a_{i}\a_{i+1},\cdots,\a_{n+1}}(a_1,\cdots, a_{i-1}, a_i*_{\a_i,\a_{i+1}}a_{i+1}, a_{i+2},\cdots, a_{n+1})\\
&&+(-1)^{n+1}f_{\a_1,\cdots, \a_{n}}(a_1,\cdots, a_n)\c_{\a_1\cdots\a_{n},\a_{n+1}}a_{n+1},\label{4.1}
\end{eqnarray}
for $f\in C_{\Omega}^{n}(A, M)$.

We use the following conventions that correspond to the composition defined in \eqref{4.1}. Suppose that $a*_{\a,\b}m=\sum_{i}f_i\o_H m_i\in H^{\o 2}\o_H M$, $m*_{\a,\b}a=\sum_{j}g_{j}\o_H n_{j}\in H^{\o 2}\o_H M$ for all $a\in A$ and $m\in M$. Then for all $f\in H^{\o n}$, we set
\begin{eqnarray*}
&&a*_{\a,\b}(f\o m)=\sum_{i}(1\o f)(\hbox{id}\o \D^{(n-1)})(f_i)\o_H m_i\in H^{\o (n+1)}\o_H M\\
\hbox{and}\\
&&(f\o m)*_{\a,\b}a=\sum_{j}(f\o 1)(\D^{(n-1)}\o \hbox{id})(g_j)\o_H n_j\in H^{\o (n+1)}\o_H M,
\end{eqnarray*}
where $\D^{(n-1)}=(\hbox{id}\o\cdots\o \hbox{id}\o\D)\cdots(\hbox{id}\o\D)\D: H\rightarrow H^{\o n}$ is the iterated comultiplication for $n>1$, and $\D^{(0)}:=\hbox{id}$.
For $g\in H^{\o 2}$ and $f=\{f_{\a_1,\cdots, \a_n}\}_{\a_1,\cdots, \a_n\in \Omega}\in C^{n}_{\Omega}(A, M)$, we set
\begin{eqnarray*}
&&f_{\a_1,\cdots,\a_{n}}(a_1\o\cdots\o a_{i-1}\o (g\o_H a_i)\o a_{i+1}\o\cdots\o a_n)\\
&=&((\hbox{id}^{\o(i-1)}\o g\D\o \hbox{id}^{\o(n-i)})\o_H \hbox{id})f_{\a_1,\cdots,\a_{n}}(a_1\o\cdots\o a_n)\in H^{\o(n+1)}\o_H M.
\end{eqnarray*}

Similar to the standard Hochschild cohomology, one can show that $d_{\Omega}^{n}\circ d_{\Omega}^{n-1}=0$. So $(C^{n}_{\Omega}(A, M),d_{\Omega}^{n})$ is a cochain complex. Let $Z_{\Omega}^{n}(A, M)=\{f\in C_{\Omega}^{n}(A, M)|d_{\Omega}^{n}f=0\}$ be the space of $n$-cocycles and $B_{\Omega}^{n}(A, M)=\{d_{\Omega}^{n-1}g|g\in C_{\Omega}^{n-1}(A, M)\}$ the space of $n$-coboundaries. The quotient space
$$
H_{\Omega}^{n}(A, M)=\frac{Z_{\Omega}^{n}(A, M)}{B_{\Omega}^{n}(A, M)}
$$
is called \textbf{the cohomology of the $\Omega$-associative $H$-pseudoalgebra} $A$ with coefficients in $M$.

In particular, by setting $n=1$ and $n=2$ in condition \eqref{4.1} respectively, we have
\begin{eqnarray}
&&(d_{\Omega}^{1}f)_{\a_1,\a_2}(a_1,a_2)=a_1\c_{\a_1,\a_2}f_{\a_2}(a_2)-f_{\a_1\a_2}(a_1*_{\a_1,\a_2}a_2)+f_{\a_1}(a_1)\c_{\a_1,\a_2}a_2\label{4.2}
\end{eqnarray}
and
\begin{eqnarray}
\nonumber&&(d_{\Omega}^{2}f)_{\a_1,\a_2,\a_3}(a_1,a_2,a_3)=a_1\c_{\a_1,\a_2\a_3}f_{\a_2,\a_3}(a_2,a_3)-f_{\a_1,\a_2}(a_1,a_2)\c_{\a_1\a_2,\a_3}a_3\\
&&\quad\quad\quad\quad\quad\quad\quad\quad\quad\quad-f_{\a_1\a_2,\a_3}(a_1*_{\a_1,\a_2}a_2,a_3)+f_{\a_1,\a_2\a_3}(a_1,a_2*_{\a_2,\a_3}a_3).\label{4.3}
\end{eqnarray}
\\

{\bf 4.2 The cohomology of pseudo-$\mathcal{O}$-operator family}\\

In this subsection, we use a pseudo-$\mathcal{O}$-operator family to induce new $\Omega$-associative $H$-pseudoalgebra structures and their bimodules, and subsequently define the cohomology of pseudo-$\mathcal{O}$-operator families.\\

{\bf Proposition 4.2} Let $(A, *)$ be an associative $H$-pseudoalgebra, $M$ an $A$-bimodule and $\{T_\a: M\rightarrow A\}_{\a\in \Omega}$ a pseudo-$\mathcal{O}$-operator family. Then $(M, \{\widehat{*}_{\a,\b}\}_{\a,\b\in \Omega})$ is an $\Omega$-associative $H$-pseudoalgebra, where
\begin{eqnarray*}
&&u\widehat{*}_{\a,\b}v=T_{\a}(u)*v+u*T_{\b}(v),
\end{eqnarray*}
for all $u, v\in M$ and $\a,\b\in \Omega$.

{\bf Proof.} It is easy to prove that $\{\widehat{*}_{\a,\b}\}_{\a,\b\in \Omega}$ are $H^{\o 2}$-linear maps. For all $u,v,w\in M$ and $\a,\b,\g\in \Omega$, we have
\begin{eqnarray*}
&&(u\widehat{*}_{\a,\b}v)\widehat{*}_{\a\b,\g}w\\
&=&(T_{\a}(u)*v+u*T_{\b}(v))\widehat{*}_{\a\b,\g}w\\
&=&T_{\a\b}(T_{\a}(u)*v+u*T_{\b}(v))*w+(T_{\a}(u)*v+u*T_{\b}(v))*T_{\g}(w)\\
&=&(T_{\a}(u)*T_{\b}(v))*w+(T_{\a}(u)*v)*T_{\g}(w)+(u*T_{\b}(v))*T_{\g}(w)\\
&=&T_{\a}(u)*(T_{\b}(v)*w+v*T_{\g}(w))+u*T_{\b\g}(T_{\b}(v)*w+v*T_{\g}(w))\\
&=&u\widehat{*}_{\a,\b\g}(T_{\b}(v)*w+v*T_{\g}(w))\\
&=&u\widehat{*}_{\a,\b\g}(v\widehat{*}_{\b,\g}w).
\end{eqnarray*}
Hence $(M, \{\widehat{*}_{\a,\b}\}_{\a,\b\in \Omega})$ is an $\Omega$-associative $H$-pseudoalgebra.
$\hfill \square$
\\

{\bf Theorem 4.3} Let $(A, *)$ be an associative $H$-pseudoalgebra, $M$ an $A$-bimodule and $\{T_\a: M\rightarrow A\}_{\a\in \Omega}$ a pseudo-$\mathcal{O}$-operator family. Define a collection of $H$-bilinear maps

\begin{align*}
\begin{Bmatrix}
\triangleright_{\a,\b}: M\o A\rightarrow A, u\triangleright_{\a,\b}a=T_\a(u)*a-T_{\a\b}(u*a)\\
\triangleleft_{\a,\b}: A\o M\rightarrow A, a\triangleleft_{\a,\b}u=a*T_\b(u)-T_{\a\b}(a*u)
\end{Bmatrix}_{\a,\b\in \Omega}.
\end{align*}
Then $(A, \{\triangleright_{\a,\b},\triangleleft_{\a,\b}\}_{\a,\b\in \Omega})$ is a bimodule over $(M, \{\widehat{*}_{\a,\b}\}_{\a,\b\in \Omega})$ defined in Proposition 4.2.

{\bf Proof.} For all $u,v\in M$, $a\in A$ and $\a,\b,\g\in \Omega$, using equation \eqref{2.1}, we have
\begin{eqnarray*}
&&(u\widehat{*}_{\a,\b}v)\triangleright_{\a\b,\g}a-u\triangleright_{\a,\b\g}(v\triangleright_{\b,\g}a)\\
&=&(T_{\a}(u)*v+u*T_{\b}(v))\triangleright_{\a\b,\g}a-u\triangleright_{\a,\b\g}(T_\b(v)*a-T_{\b\g}(v*a))\\
&=&T_{\a\b}(T_{\a}(u)*v+u*T_{\b}(v))*a-T_{\a\b\g}((T_{\a}(u)*v)*a)-T_{\a\b\g}((u*T_{\b}(v))*a)\\
&&-T_{\a}(u)*(T_{\b}(v)*a)+T_{\a}(u)*T_{\b\g}(v*a)+T_{\a\b\g}(u*(T_{\b}(v)*a))-T_{\a\b\g}(u*T_{\b\g}(v*a))\\
&=&-T_{\a\b\g}((T_{\a}(u)*v)*a)+T_{\a}(u)*T_{\b\g}(v*a)-T_{\a\b\g}(u*T_{\b\g}(v*a))\\
&=&0,
\end{eqnarray*}

\begin{eqnarray*}
&&(u\triangleright_{\a,\b}a)\triangleleft_{\a\b,\g}v-u\triangleright_{\a,\b\g}(a\triangleleft_{\b,\g}v)\\
&=&(T_\a(u)*a-T_{\a\b}(u*a))\triangleleft_{\a\b,\g}v-u\triangleright_{\a,\b\g}(a*T_\g(v)-T_{\b\g}(a*v))\\
&=&(T_\a(u)*a)*T_\g(v)-T_{\a\b}(u*a)*T_{\g}(v)-T_{\a\b\g}((T_{\a}(u)*a)*v)+T_{\a\b\g}(T_{\a\b}(u*a)*v)\\
&&-T_{\a}(u)*(a*T_{\g}(v))+T_{\a}(u)*T_{\b\g}(a*v)+T_{\a\b\g}(u*(a*T_{\g}(v)))-T_{\a\b\g}(u*T_{\b\g}(a*v))\\
&=&0,
\end{eqnarray*}
and
\begin{eqnarray*}
&&(a\triangleleft_{\a,\b}u)\triangleleft_{\a\b,\g}v-a\triangleleft_{\a,\b\g}(u\widehat{*}_{\b,\g}v)\\
&=&(a*T_\b(u)-T_{\a\b}(a*u))\triangleleft_{\a\b,\g}v-a\triangleleft_{\a,\b\g}(T_{\b}(u)*v+u*T_{\g}(v))\\
&=&(a*T_\b(u))*T_{\g}(v)-T_{\a\b}(a*u)*T_{\g}(v)-T_{\a\b\g}((a*T_{\b}(u))*v)+T_{\a\b\g}(T_{\a\b}(a*u)*v)\\
&&-a*T_{\b\g}(T_{\b}(u)*v)-a*T_{\b\g}(u*T_{\g}(v))+T_{\a\b\g}(a*(T_{\b}(u)*v))+T_{\a\b\g}(a*(u*T_{\g}(v)))\\
&=&a*(T_\b(u)*T_{\g}(v))-a*T_{\b\g}(T_{\b}(u)*v)-a*T_{\b\g}(u*T_{\g}(v))\\
&&-T_{\a\b}(a*u)*T_{\g}(v)+T_{\a\b\g}(T_{\a\b}(a*u)*v)+T_{\a\b\g}((a*u)*T_{\g}(v))\\
&=&0.
\end{eqnarray*}
So we complete the proof.
$\hfill \square$
\\

It follows that one can consider the cochain complex of the $\Omega$-associative $H$-pseudoalgebra $(M, \{\widehat{*}_{\a,\b}\}_{\a,\b\in \Omega})$ with
coefficients in the bimodule given in the above theorem. More precisely, we define

\begin{eqnarray*}
&&C_{\Omega}^{0}(M, A)=k\o_H A,\\
&&C_{\Omega}^{n\geq1}(M, A)=\{f=\{f_{\a_1,\cdots,\a_n}\}_{\a_1,\cdots,\a_n\in \Omega}\in \hbox{Hom}_{H^{\o n}}(M^{\o n}, H^{\o n}\o_H A)\}
\end{eqnarray*}
and the differential $\delta^{n}: C_{\Omega}^{n}(M, A)\rightarrow C_{\Omega}^{n+1}(M, A)$ is given by
\begin{eqnarray*}
&&\delta^{0}(1\o a)_{\a}(u)=\mu_{1,\v}(T_\a(u)*a-T_\a(u*a))-\mu_{\v,1}(a*T_\a(u)-T_\a(a*u)),
\end{eqnarray*}
where $\mu_{1,\v}(h\o g\o_H a)=\v(g)h\c a$ and $\mu_{\v,1}(h\o g\o_H a)=\v(h)g\c a$, for all $h,g\in H$ and $a\in A$.

\begin{eqnarray*}
&&(\delta^{n\geq1}f)_{\a_1,\cdots,\a_{n+1}}(u_1,\cdots, u_{n+1})\\
&=&T_{\a_1}(u_1)*f_{\a_2,\cdots, \a_{n+1}}(u_2,\cdots, u_{n+1})-T_{\a_1\cdots\a_{n+1}}(u_1*f_{\a_2,\cdots,\a_{n+1}}(u_2,\cdots,u_{n+1}))\\
&&+\sum_{i=1}^{n}(-1)^{i}f_{\a_1,\cdots, \a_{i}\a_{i+1},\cdots,\a_{n+1}}(u_1,\cdots, u_{i-1}, T_{\a_i}(u_i)*u_{i+1}+u_{i}*T_{\a_{i+1}}(u_{i+1}),\\
&&u_{i+2},\cdots, u_{n+1})+(-1)^{n+1}f_{\a_1,\cdots, \a_{n}}(u_1,\cdots, u_n)*T_{\a_{n+1}}(u_{n+1})\\
&&-(-1)^{n+1}T_{\a_1\cdots\a_{n+1}}(f_{\a_1,\cdots, \a_{n}}(u_1,\cdots, u_n)*u_{n+1}),
\end{eqnarray*}
for $f\in C_{\Omega}^{n}(M, A)$.

The corresponding cohomology groups are called \textbf{the cohomology of the pseudo-$\mathcal{O}$-operator family} $\{T_\a\}_{\a\in \Omega}$.

\section*{5. Deformation theory of $\Omega$-associative $H$-pseudoalgebras}
\def\theequation{5. \arabic{equation}}
\setcounter{equation} {0} \hskip\parindent

In this section, we study the formal deformations of $\Omega$-associative $H$-pseudoalgebras and characterize them in terms of lower degree cohomology
groups of $\Omega$-associative $H$-pseudoalgebras.\\

{\bf Definition 5.1} Let $(A,(*_{\a,\b})_{\a,\b\in \Omega})$ be an $\Omega$-associative $H$-pseudoalgebra. A \textbf{formal deformation} of $A$ consists of a formal sum
\begin{eqnarray*}
&&\mathfrak{T}^{t}_{\a,\b}=\sum_{i=0}^{\infty}T^{(i)}_{\a,\b}t^{i}=T^{(0)}_{\a,\b}+T^{(1)}_{\a,\b}t+T^{(2)}_{\a,\b}t^{2}+\cdots,
\end{eqnarray*}
where each $T^{(i)}=(T^{(i)}_{\a,\b})_{\a,\b\in \Omega}: A\o A\rightarrow H^{\o2}\o_H A$ is an $H^{\o 2}$-linear map and $T^{(0)}_{\a,\b}(x,y)=x*_{\a,\b}y$, such that $\mathfrak{T}^{t}_{\a,\b}$ endows $A$ with an $\Omega$-associative $H$-pseudoalgebra structure. That is,
\begin{eqnarray}
&&\mathfrak{T}^{t}_{\a,\b\g}(x, \mathfrak{T}^{t}_{\b,\g}(y,z))=\mathfrak{T}^{t}_{\a\b,\g}(\mathfrak{T}^{t}_{\a,\b}(x, y),z), \label{5.1}
\end{eqnarray}
for all $x, y,z\in A$ and $\a,\b,\g\in \Omega$. This system of equations is called \textbf{the deformation equations}.

Comparing the coefficients of $t^{n}$ in \eqref{5.1}, we obtain that $(T^{(i)}_{\a,\b})_{i\geq0,\a,\b\in \Omega}$ satisfies
\begin{eqnarray*}
&&\sum_{i=0}^{s}T^{(i)}_{\a,\b\g}(x,T^{(s-i)}_{\b,\g}(y,z))=\sum_{i=0}^{s}T^{(i)}_{\a\b,\g}(T^{(s-i)}_{\a,\b}(x,y),z).
\end{eqnarray*}
We denote $T^{(i)}\hat{\circ}T^{(j)}(x,y,z)_{\a,\b,\g}=T^{(i)}_{\a,\b\g}(x,T^{(j)}_{\b,\g}(y,z))-T^{(i)}_{\a\b,\g}(T^{(j)}_{\a,\b}(x,y),z)$, then the above equation can be written as
\begin{eqnarray*}
&&\sum_{i=0}^{s}T^{(i)}\hat{\circ}T^{(s-i)}=0, \quad s=0,1,2,\cdots.
\end{eqnarray*}

Obviously, for $s=0$, this corresponds to the associativity of $(A,(*_{\a,\b})_{\a,\b\in \Omega})$.

For $s=1$, the equation $T^{(0)}\hat{\circ} T^{(1)}+T^{(1)}\hat{\circ} T^{(0)}=0$ is equivalent to
\begin{eqnarray}
&&x*_{\a,\b\g}T^{(1)}_{\b,\g}(y,z)+T^{(1)}_{\a,\b\g}(x,y*_{\b,\g}z)=T^{(1)}_{\a,\b}(x,y)*_{\a\b,\g}z+T^{(1)}_{\a\b,\g}(x*_{\a,\b}y,z).\label{5.2}
\end{eqnarray}

For $s=2$, we have $T^{(0)}\hat{\circ} T^{(2)}+T^{(1)}\hat{\circ} T^{(1)}+T^{(2)}\hat{\circ} T^{(0)}=0$. This implies a key relation for the coboundary operator
$$
d_{\Omega}^{2}T^{(2)}=T^{(0)}\hat{\circ}T^{(2)}+T^{(2)}\hat{\circ}T^{(0)}=-T^{(1)}\hat{\circ}T^{(1)}.
$$

$\cdots$

For a general $s=n$, we obtain the general deformation equation

$d_{\Omega}^{2}T^{(n)}=T^{(0)}\hat{\circ}T^{(n)}+T^{(n)}\hat{\circ}T^{(0)}=-(T^{(1)}\hat{\circ}T^{(n-1)}+T^{(2)}\hat{\circ}T^{(n-2)}+\cdots+T^{(n-1)}\hat{\circ}T^{(1)})$.

 Note that $(d_{\Omega}^{2}T^{(1)})_{\a,\b,\g}=0$ by conditions \eqref{4.3} and \eqref{5.2}. This implies that $T^{(1)}=(T^{(1)}_{\a,\b})_{\a,\b\in \Omega}$ is a 2-cocycle in the cohomology of $(A, (*_{\a,\b})_{\a,\b\in \Omega})$. We call $(T^{(1)}_{\a,\b})_{\a,\b\in \Omega}$ \textbf{the infinitesimal deformation} of $A$.\\

{\bf Definition 5.2} Let $(A,(*_{\a,\b})_{\a,\b\in \Omega})$ be an $\Omega$-associative $H$-pseudoalgebra. Two deformations $\mathfrak{T}^{t}_{\a,\b}=\sum_{i=0}^{\infty}T^{(i)}_{\a,\b}t^{i}$ and $\mathfrak{\widetilde{T}}^{t}_{\a,\b}=\sum_{i=0}^{\infty}\widetilde{T}^{(i)}_{\a,\b}t^{i}$ of $A$
are said to be \textbf{equivalent} if there exists a family $(\Phi^{t}_{\omega})_{\omega\in \Omega}: A\rightarrow A[[t]]$ of $H$-linear maps of the form $\Phi^{t}_{\omega}=\sum_{i=0}^{\infty}\phi_{\omega}^{(i)}t^{i}=\phi^{(0)}_{\omega}+\phi^{(1)}_{\omega}t^{1}+\phi^{(2)}_{\omega}t^{2}+\cdots$, where each $\phi^{(i)}_{\omega}: A\rightarrow A$ is an $H$-linear map and $\phi^{(0)}_{\omega}=\hbox{id}$, such that
\begin{eqnarray*}
&&\Phi^{t}_{\a\b}\circ\mathfrak{T}^{t}_{\a,\b}(x,y)=\mathfrak{\widetilde{T}}^{t}_{\a,\b}(\Phi^{t}_{\a}(x),\Phi^{t}_{\b}(y)).
\end{eqnarray*}
More precisely,
\begin{eqnarray}
&&\sum_{i=0}^{\infty}\phi^{(i)}_{\a\b}(\sum_{j=0}^{\infty}T^{(j)}_{\a,\b}(x, y)t^{j})t^{i}=\sum_{i=0}^{\infty}\widetilde{T}^{(i)}_{\a,\b}(\sum_{j=0}^{\infty}\phi^{(j)}_{\a}(x)t^{j},\sum_{k=0}^{\infty}\phi^{(k)}_{\b}(y)t^{k})t^{i}.\label{5.3}
\end{eqnarray}

By equating coefficients of $t$ in both sides of \eqref{5.3}, we obtain
\begin{eqnarray*}
&&T^{(1)}_{\a,\b}(x,y)+\phi^{(1)}_{\a\b}(x*_{\a,\b}y)=\widetilde{T}^{(1)}_{\a,\b}(x,y)+\phi^{(1)}_{\a}(x)*_{\a,\b}y+x*_{\a,\b}\phi^{(1)}_{\b}(y),
\end{eqnarray*}
that is,
\begin{eqnarray*}
(T^{(1)}_{\a,\b}-\widetilde{T}^{(1)}_{\a,\b})(x,y)&=&\phi^{(1)}_{\a}(x)*_{\a,\b}y+x*_{\a,\b}\phi^{(1)}_{\b}(y)-\phi^{(1)}_{\a\b}(x*_{\a,\b}y)\\
&=&(d_{\Omega}^{1}\phi^{(1)})_{\a,\b}(x,y)\quad(\hbox{by Eq. \eqref{4.2}}).
\end{eqnarray*}
As a summary, we get the following theorem.\\

{\bf Theorem 5.3} Let $\mathfrak{T}^{t}_{\a,\b}$ and $\mathfrak{\widetilde{T}}^{t}_{\a,\b}$ be two equivalent formal deformation of the $\Omega$-associative $H$-pseudoalgebra $(A,(*_{\a,\b})_{\a,\b\in \Omega})$. Then there is a one-to-one correspondence between the elements of $H_{\Omega}^{2}(A,A)$ and the infinitesimal deformations of $A$.
\\

An $\Omega$-associative $H$-pseudoalgebra $(A,(*_{\a,\b})_{\a,\b\in \Omega})$ is called \textbf{rigid}, if every formal deformation $\mathfrak{T}^{t}_{\a,\b}$ of $A$ is equivalent to $\mathfrak{\widetilde{T}}^{t}_{\a,\b}=T^{(0)}_{\a,\b}$.
Now we characterize the rigidity of the $\Omega$-associative $H$-pseudoalgebra in terms of its second cohomology group.\\

{\bf Theorem 5.4} Let $(A,(*_{\a,\b})_{\a,\b\in \Omega})$ be an $\Omega$-associative $H$-pseudoalgebra and $H_{\Omega}^{2}(A,A)=0$. Then $(A,(*_{\a,\b})_{\a,\b\in \Omega})$ is rigid.

{\bf Proof.} Let $\mathfrak{T}^{t}_{\a,\b}$ be the deformation of $(A,(*_{\a,\b})_{\a,\b\in \Omega})$. Suppose that $\mathfrak{T}^{t}_{\a,\b}=T^{(0)}_{\a,\b}+\sum_{i\geq n}T^{(i)}_{\a,\b}t^{i}$, then
\begin{eqnarray*}
&&d_{\Omega}^{2}T^{(n)}=-(T^{(1)}\hat{\circ}T^{(n-1)}+T^{(2)}\hat{\circ}T^{(n-2)}+\cdots+T^{(n-1)}\hat{\circ}T^{(1)})=0.
\end{eqnarray*}
Hence $T^{(n)}=\{T^{(n)}_{\a,\b}\}_{\a,\b\in \Omega}\in Z^{2}(A,A)=B^{2}(A,A)$. It follows that there exists $\phi=\{\phi_{\omega}\}_{\omega\in \Omega}\in C_{\Omega}^{1}(A,A)$ such that $T^{(n)}_{\a,\b}=(d_{\Omega}^{1}\phi)_{\a,\b}$.

Let $(\Phi^{t}_{\omega})_{\omega\in \Omega}=\hbox{id}-\phi_{\omega} t^{n}:A\rightarrow A$. Note that each $\Phi^{t}_{\omega}$ is an isomorphism of $H$-modules with
\begin{eqnarray*}
&&\Phi^{t}_{\omega}\circ\sum_{i\geq0}\phi^{i}_{\omega}t^{in}=\sum_{i\geq0}\phi^{i}_{\omega}t^{in}\circ\Phi^{t}_{\omega}=\hbox{id}_{A}.
\end{eqnarray*}
 Define
 $$
 \mathfrak{\widetilde{T}}^{t}_{\a,\b}(x,y)=(\Phi^{t}_{\a\b})^{-1}\mathfrak{T}^{t}_{\a,\b}(\Phi^{t}_{\a}(x),\Phi^{t}_{\b}(y)).
 $$
 It is immediate that $\mathfrak{\widetilde{T}}^{t}_{\a,\b}$ is a formal deformation of $A$ and $\mathfrak{T}^{t}_{\a,\b}$ is equivalent to $\mathfrak{\widetilde{T}}^{t}_{\a,\b}$. Suppose that $\mathfrak{\widetilde{T}}^{t}_{\a,\b}=\widetilde{T}^{(0)}_{\a,\b}+\sum_{i>0}\widetilde{T}^{(i)}_{\a,\b}t^{i}$. Then

\begin{eqnarray*}
&&(\hbox{id}-\phi_{\a\b} t^{n})(\widetilde{T}^{(0)}_{\a,\b}(x,y)+\sum_{i>0}\widetilde{T}^{(i)}_{\a,\b}(x,y)t^{i})\\
&=&(T^{(0)}_{\a,\b}+\sum_{i\geq n}T^{(i)}_{\a,\b}t^{i})(x-\phi_{\a}(x)t^{n},y-\phi_{\b}(y)t^{n}),
\end{eqnarray*}
that is,
\begin{eqnarray*}
&&x*_{\a,\b}y+\sum_{i>0}\widetilde{T}^{(i)}_{\a,\b}(x,y)t^{i}-\phi_{\a\b}(x*_{\a,\b}y+\sum_{i>0}\widetilde{T}^{(i)}_{\a,\b}(x,y)t^{i})t^{n}\\
&=&x*_{\a,\b}y-(x*_{\a,\b}\phi_{\b}(y)+\phi_{\a}(x)*_{\a,\b}y)t^{n}+(\phi_{\a}(x)*_{\a,\b}\phi_{\b}(y))t^{2n}+\sum_{i\geq n}T^{(i)}_{\a,\b}(x,y)t^{i}\\
&&-\sum_{i\geq n}(T^{(i)}_{\a,\b}(\phi_{\a}(x),y)+T^{(i)}_{\a,\b}(x,\phi_{\b}(y)))t^{i+n}+\sum_{i\geq n}T^{(i)}_{\a,\b}(\phi_{\a}(x),\phi_{\b}(y))t^{i+2n}.
\end{eqnarray*}
Then we have $\widetilde{T}^{(1)}_{\a,\b}=\widetilde{T}^{(2)}_{\a,\b}=\cdots=\widetilde{T}^{(n-1)}_{\a,\b}=0$ and the coefficients of $t^{n}$ are identical, i.e.,
\begin{eqnarray*}
&&\widetilde{T}^{(n)}_{\a,\b}(x,y)-\phi_{\a\b}(x*_{\a,\b}y)=-(\phi_{\a}(x)*_{\a,\b}y+x*_{\a,\b}\phi_{\b}(y))+T^{(n)}_{\a,\b}(x,y).
\end{eqnarray*}
It follows that $\widetilde{T}^{(n)}_{\a,\b}(x,y)=T^{(n)}_{\a,\b}(x,y)-(d_{\Omega}^{1}\phi)_{\a,\b}(x,y)=0$ and $\mathfrak{\widetilde{T}}^{t}_{\a,\b}=\widetilde{T}^{(0)}_{\a,\b}+\sum_{i\geq n}\widetilde{T}^{(i)}_{\a,\b}t^{i}$. By induction, this procedure ends with $\mathfrak{\widetilde{T}}^{t}_{\a,\b}=\widetilde{T}^{(0)}_{\a,\b}$. So $\mathfrak{T}^{t}_{\a,\b}$ is equivalent to $\widetilde{T}^{(0)}_{\a,\b}$, as desired.
$\hfill \square$

\section*{6. $\Omega$-Poisson $H$-pseudoalgebras}
\def\theequation{6. \arabic{equation}}
\setcounter{equation} {0} \hskip\parindent

In this section, we show that the first-order element of a deformation, as defined in the previous section, induces an $\Omega$-Poisson $H$-pseudoalgebra. We begin by defining this concept.\\

{\bf Definition 6.1} An \textbf{$\Omega$-Poisson $H$-pseudoalgebra} is a triple $(P, ([\cdot*_{\a,\b}\cdot])_{\a,\b\in \Omega}, (*_{\a,\b})_{\a,\b\in \Omega})$, such that $(P, ([\cdot*_{\a,\b}\cdot])_{\a,\b\in \Omega})$ is an $\Omega$-Lie $H$-pseudoalgebra, $(P, (*_{\a,\b})_{\a,\b\in \Omega})$ is a commutative $\Omega$-associative $H$-pseudoalgebra and they satisfy the compatible condition
\begin{eqnarray}
&&[x*_{\a,\b\g}(y*_{\b,\g}z)]=[x*_{\a,\b}y]*_{\a\b,\g}z+((12)\o_H id)y*_{\b,\a\g}[x*_{\a,\g}z], \label{6.1}
\end{eqnarray}
for all $x, y, z\in P$ and $\a,\b,\g\in \Omega$.
\\

{\bf Remark 6.2} When $\Omega=\{0\}$, an $\Omega$-Poisson $H$-pseudoalgebra recovers the notion of a Poisson $H$-pseudoalgebra. Furthermore, if $H=k$, we obtain an ordinary Poisson algebra.\\

As per the definition, any $\Omega$-Poisson $H$-pseudoalgebra has an underlying commutative $\Omega$-associative $H$-pseudoalgebra structure. Conversely, any commutative $\Omega$-associative $H$-pseudoalgebra $(A, (*_{\a,\b})_{\a,\b\in \Omega})$ can be equipped with the structure of an $\Omega$-Poisson $H$-pseudo-\\
algebra by defining $[\cdot *_{\alpha,\beta} \cdot]=0$ for all $\alpha, \beta \in \Omega$. Moreover, we have the following conclusion.\\

{\bf Proposition 6.3} Let $(A, (*_{\a,\b})_{\a,\b\in \Omega})$ be a commutative $\Omega$-associative $H$-pseudoalgebra. Define $[x*_{\a,\b}y]=x*_{\a,\b}y-(\sigma\o_H \hbox{id})(y*_{\b,\a}x)$ for all $x, y\in A$, then $(A, (*_{\a,\b})_{\a,\b\in \Omega}, ([\cdot*_{\a,\b}\cdot])_{\a,\b\in \Omega})$ is an $\Omega$-Poisson $H$-pseudoalgebra.

{\bf Proof.} By Corollary 3.20, $(A, [\cdot*_{\a,\b}\cdot])$ is an $\Omega$-Lie $H$-pseudoalgebra. It is sufficient to verify the compatibility condition \eqref{6.1}. For all $x, y, z\in A$ and $\a, \b, \g\in \Omega$, we have
\begin{eqnarray*}
&&[x*_{\a,\b}y]*_{\a\b,\g}z+((12)\o_H \hbox{id})y*_{\b,\a\g}[x*_{\a,\g}z]\\
&=&(x*_{\a,\b}y-(\sigma\o_H \hbox{id})(y*_{\b,\a}x))*_{\a\b,\g}z+((12)\o_H \hbox{id})y*_{\b,\a\g}(x*_{\a,\g}z)\\
&&-((12)\o_H \hbox{id})(y*_{\b,\a\g}(\sigma\o_H \hbox{id})(z*_{\g,\a}x))\\
&=&(x*_{\a,\b}y)*_{\a\b,\g}z-((12)\o_H \hbox{id})(y*_{\b,\a}x)*_{\a\b,\g}z+((12)\o_H \hbox{id})(y*_{\b,\a\g}(x*z)\\
&&-((23)\o_H \hbox{id})y*_{\b,\a\g}(z*_{\g,\a}x))\\
&=&(x*_{\a,\b}y)*_{\a\b,\g}z-((12)\o_H \hbox{id})(y*_{\b,\a}x)*_{\a\b,\g}z+((12)\o_H \hbox{id})y*_{\b,\a\g}(x*_{\a,\g}z)\\
&&-((123)\o_H \hbox{id})y*_{\b,\a\g}(z*_{\g,\a}x))\\
&=&x*_{\a,\b\g}(y*_{\b,\g}z)-((123)\o_H \hbox{id})(y*_{\b,\g}z)*_{\b\g,\a}x\\
&=&[x*_{\a,\b\g}(y*_{\b,\g}z)],
\end{eqnarray*}
which completes the proof.
$\hfill \square$
\\

To prove the main result of this section, we need the following lemmas.\\

{\bf Lemma 6.4} Let $(A, (*_{\a,\b})_{\a,\b\in \Omega})$ be a commutative $\Omega$-associative $H$-pseudoalgebra, and let $\phi$ be a skew-commutative 2-cochain (i.e.,
 $\phi_{\a,\b}(x,y)=-(\sigma\o_H \hbox{id})\phi_{\b,\a}(y,x), ~\forall \a,\b\in \Omega, x,y\in A$).
 Assume that $(d_{\Omega}^{2}\phi)_{\a,\b,\g}=0$. Then the following identity holds for all $x, y,z\in A$ and $\a,\b,\g\in \Omega$:
\begin{eqnarray*}
&&\phi_{\a,\b\g}(x,y*_{\b,\g}z)=\phi_{\a,\b}(x, y)*_{\a\b,\g}z+((12)\o_H \hbox{id})(y*_{\b,\a\g}\phi_{\a,\g}(x,z)).
\end{eqnarray*}

{\bf Proof.} Applying the skew-commutativity of $\phi$ and the commutativity of $A$ yields the following identities:
\begin{eqnarray}
&&\phi_{\b,\a\g}(y,x*_{\a,\g}z)=-((123)\o_H \hbox{id})\phi_{\a\g,\b}(x*_{\a,\g}z,y),\label{L1}\\
&&x*_{\a,\b\g}\phi_{\b,\g}(y,z)=-((23)\o_H \hbox{id})x*_{\a,\b\g}\phi_{\g,\b}(z,y),\label{L2}\\
&&\phi_{\a,\g}(x,z)*_{\a\g,\b}y=((132)\o_H \hbox{id})(y*_{\b,\a\g}\phi_{\a,\g}(x,z)),\label{L3}\\
&&\phi_{\a,\b}(x,y)*_{\a\b,\g}z=-((12)\o_H \hbox{id})\phi_{\b,\a}(y,x)*_{\a\b,\g}z\label{L4},
\end{eqnarray}
for all $x, y, z\in A$.
The condition $(d_2\phi)_{\a,\b,\g}=0$ implies that
\begin{eqnarray}
&&\phi_{\a,\b\g}(x,y*_{\b,\g}z)=\phi_{\a\b,\g}(x*_{\a,\b}y,z)-x*_{\a,\b\g}\phi_{\b,\g}(y,z)+\phi_{\a,\b}(x,y)*_{\a\b,\g}z,~~\label{L5}
\end{eqnarray}
Replacing $(y,z)$ for $(z,y)$ and $(\b,\g)$ for $(\g,\b)$ in \eqref{L5}, we have
\begin{eqnarray}
\nonumber&&((23)\o_H \hbox{id})\phi_{\a,\b\g}(x,z*_{\g,\b}y)\\
&&=((23)\o_H \hbox{id})(\phi_{\a\g,\b}(x*_{\a,\g}z,y)-x*_{\a,\b\g}\phi_{\g,\b}(z,y)+\phi_{\a,\g}(x,z)*_{\a\g,\b}y).~\label{L6}
\end{eqnarray}
Replacing $(x,y)$ for $(y,x)$ and $(\a,\b)$ for $(\b,\a)$ in \eqref{L5}, we obtain
\begin{eqnarray}
\nonumber&&((12)\o_H \hbox{id})(\phi_{\a\b,\g}(y*_{\b,\a}x,z)-\phi_{\b,\a\g}(y,x*_{\a,\g}z))\\
&&=((12)\o_H \hbox{id})(y*_{\b,\a\g}\phi_{\a,\g}(x,z)-\phi_{\b,\a}(y,x)*_{\a\b,\g}z).~\label{L7}
\end{eqnarray}
Adding equations \eqref{L5}-\eqref{L7}, the left-hand side becomes
\begin{eqnarray*}
&&\phi_{\a,\b\g}(x,y*_{\b,\g}z)+((23)\o_H \hbox{id})\phi_{\a,\b\g}(x,z*_{\g,\b}y)+((12)\o_H \hbox{id})\phi_{\a\b,\g}(y*_{\b,\a}x,z)\\
&&-((12)\o_H \hbox{id})\phi_{\b,\a\g}(y,x*_{\a,\g}z)\\
&=&\phi_{\a,\b\g}(x,y*_{\b,\g}z)+\phi_{\a,\b\g}(x,y*_{\b,\g}z)+\phi_{\a\b,\g}(x*_{\a,\b}y,z)-((12)\o_H \hbox{id})\phi_{\b,\a\g}(y,x*_{\a,\g}z)\\
&=&2\phi_{\a,\b\g}(x,y*_{\b,\g}z)+\phi_{\a\b,\g}(x*_{\a,\b}y,z)\\
&&+((12)(123)\o_H \hbox{id})\phi_{\a\g,\b}(x*_{\a,\g}z,y) \quad(\hbox{by Eq. \eqref{L1}})\\
&=&2\phi_{\a,\b\g}(x,y*_{\b,\g}z)+\phi_{\a\b,\g}(x*_{\a,\b}y,z)+((23)\o_H \hbox{id})\phi_{\a\g,\b}(x*_{\a,\g}z,y),
\end{eqnarray*}
while the right-hand side gives
\begin{eqnarray*}
&&\phi_{\a\b,\g}(x*_{\a,\b}y,z)-x*_{\a,\b\g}\phi_{\b,\g}(y,z)+\phi_{\a,\b}(x,y)*_{\a\b,\g}z\\
&&+((23)\o_H \hbox{id})(\phi_{\a\g,\b}(x*_{\a,\g}z,y)-x*_{\a,\b\g}\phi_{\g,\b}(z,y)+\phi_{\a,\g}(x,z)*_{\a\g,\b}y)\\
&&+((12)\o_H \hbox{id})(y*_{\b,\a\g}\phi_{\a,\g}(x,z)-\phi_{\b,\a}(y,x)*_{\a\b,\g}z)\\
&=&\phi_{\a\b,\g}(x*_{\a,\b}y,z)+((23)\o_H \hbox{id})(x*_{\a,\b\g}\phi_{\g,\b}(z,y))+\phi_{\a,\b}(x,y)*_{\a\b,\g}z\\
&&+((23)\o_H \hbox{id})\phi_{\a\g,\b}(x*_{\a,\g}z,y)-((23)\o_H \hbox{id})(x*_{\a,\b\g}\phi_{\g,\b}(z,y))\\
&&+((23)(132)\o_H \hbox{id})(y*_{\b,\a\g}\phi_{\a,\g}(x,z))+((12)\o_H \hbox{id})(y*_{\b,\a\g}\phi_{\a,\g}(x,z))\\
&&+\phi_{\a,\b}(x,y)*_{\a\b,\g}z  \quad(\hbox{by Eqs. \eqref{L2}-\eqref{L4}})\\
&=&\phi_{\a\b,\g}(x*_{\a,\b}y,z)+((23)\o_H \hbox{id})\phi_{\a\g,\b}(x*_{\a,\g}z,y)+2\phi_{\a,\b}(x,y)*_{\a\b,\g}z\\
&&+2((12)\o_H \hbox{id})(y*_{\b,\a\g}\phi_{\a,\g}(x,z)).
\end{eqnarray*}
It follows that
\begin{eqnarray*}
&&2\phi_{\a,\b\g}(x,y*_{\b,\g}z)+\phi_{\a\b,\g}(x*_{\a,\b}y,z)+((23)\o_H \hbox{id})\phi_{\a\g,\b}(x*_{\a,\g}z,y)\\
&=&\phi_{\a\b,\g}(x*_{\a,\b}y,z)+((23)\o_H \hbox{id})\phi_{\a\g,\b}(x*_{\a,\g}z,y)+2\phi_{\a,\b}(x,y)*_{\a\b,\g}z\\
&&+2((12)\o_H \hbox{id})(y*_{\b,\a\g}\phi_{\a,\g}(x,z)),
\end{eqnarray*}
that is,
\begin{eqnarray*}
&&\phi_{\a,\b\g}(x,y*_{\b,\g}z)=\phi_{\a,\b}(x, y)*_{\a\b,\g}z+((12)\o_H \hbox{id})(y*_{\b,\a\g}\phi_{\a,\g}(x,z)).
\end{eqnarray*}
So we complete the proof.
$\hfill \square$
\\

{\bf Lemma 6.5} Let $(A, (*_{\a,\b})_{\a,\b\in \Omega})$ be a commutative $\Omega$-associative $H$-pseudoalgebra and $\mathfrak{T}^{t}_{\a,\b}=\sum_{i=0}^{\infty}T^{(i)}_{\a,\b}t^{i}$ the deformation of $A$. Then
\begin{eqnarray*}
&&-T^{(1)}\hat{\circ} T^{(1)}(x,y,z)_{\a,\b,\g}-((123)\o_H \hbox{id})T^{(1)}\hat{\circ} T^{(1)}(y,z,x)_{\b,\g,\a}\\
&&-((132)\o_H \hbox{id})T^{(1)}\hat{\circ} T^{(1)}(z,x,y)_{\g,\a,\b}\\
&=&T^{(2)}\hat{\circ} T^{(0)}(x,y,z)_{\a,\b,\g}+((123)\o_H \hbox{id})T^{(2)}\hat{\circ} T^{(0)}(y,z,x)_{\b,\g,\a}\\
&&+((132)\o_H \hbox{id})T^{(2)}\hat{\circ} T^{(0)}(z,x,y)_{\g,\a,\b}.
\end{eqnarray*}
for all $x, y, z\in A$ and $\a,\b,\g\in \Omega$.

{\bf Proof.} For all $x, y, z\in A$ and $\a,\b,\g\in \Omega$, we have the following identities by a direct computation:
\begin{eqnarray}
&&x*_{\a,\b\g}T^{(2)}_{\b,\g}(y,z)=((123)\o_H \hbox{id})(T^{(2)}_{\b,\g}(y,z)*_{\b\g,\a}x),\label{L8}\\
&&T^{(2)}_{\b,\g}(y,z)*_{\b\g,\a}x=((132)\o_H \hbox{id})(x*_{\a,\b\g}T^{(2)}_{\b,\g}(y,z)),\label{L9}\\
&&T^{(2)}_{\b,\a\g}(y,z*_{\g,\a}x)=((23)\o_H \hbox{id})T^{(2)}_{\b,\a\g}(y,x*_{\a,\g}z),\label{L10}\\
&&T^{(2)}_{\a\g,\b}(z*_{\g,\a}x,y)=((12)\o_H \hbox{id})T^{(2)}_{\a\g,\b}(x*_{\a,\g}z,y).\label{L11}
\end{eqnarray}
Since $T^{(0)}\hat{\circ}T^{(2)}+T^{(1)}\hat{\circ} T^{(1)}+T^{(2)}\hat{\circ} T^{(0)}=0$, we have
\begin{eqnarray}
\nonumber-T^{(1)}\hat{\circ} T^{(1)}(x, y, z)_{\a,\b,\g}&=&T^{(0)}\hat{\circ} T^{(2)}(x, y, z)_{\a,\b,\g}+T^{(2)}\hat{\circ} T^{(0)}(x, y, z)_{\a,\b,\g}\\
\nonumber&=&x*_{\a,\b\g}T^{(2)}_{\b,\g}(y,z)-T^{(2)}_{\a,\b}(x,y)*_{\a\b,\g}z+T^{(2)}_{\a,\b\g}(x, y*_{\b,\g}z)\\
&&-T^{(2)}_{\a\b,\g}(x*_{\a,\b}y,z).\label{L12}
\end{eqnarray}
For all $x, y, z\in A$ and $\a,\b,\g\in \Omega$, we compute as follows:
\begin{eqnarray*}
&&-T^{(1)}\hat{\circ} T^{(1)}(x,y,z)_{\a,\b,\g}-((123)\o_H \hbox{id})T^{(1)}\hat{\circ} T^{(1)}(y,z,x)_{\b,\g,\a}\\
&&-((132)\o_H \hbox{id})T^{(1)}\hat{\circ} T^{(1)}(z,x,y)_{\g,\a,\b}\\
&=&x*_{\a,\b\g}T^{(2)}_{\b,\g}(y,z)-T^{(2)}_{\a,\b}(x,y)*_{\a\b,\g}z+T^{(2)}_{\a,\b\g}(x, y*_{\b,\g}z)-T^{(2)}_{\a\b,\g}(x*_{\a,\b}y,z)\\
&&+((123)\o_H \hbox{id})(y*_{\b,\a\g}T^{(2)}_{\g,\a}(z,x)-T^{(2)}_{\b,\g}(y,z)*_{\b\g,\a}x+T^{(2)}_{\b,\a\g}(y, z*_{\g,\a}x)\\
&&-T^{(2)}_{\b\g,\a}(y*_{\b,\g}z,x))+((132)\o_H \hbox{id})(z*_{\g,\a\b}T^{(2)}_{\a,\b}(x,y)-T^{(2)}_{\g,\a}(z,x)*_{\g\a,\b}y\\
&&+T^{(2)}_{\g,\a\b}(z, x*_{\a,\b}y)-T^{(2)}_{\g\a,\b}(z*_{\g,\a}x,y))\quad(\hbox{by Eq. \eqref{L12}})\\
&=&x*_{\a,\b\g}T^{(2)}_{\b,\g}(y,z)-T^{(2)}_{\a,\b}(x,y)*_{\a\b,\g}z+T^{(2)}_{\a,\b\g}(x, y*_{\b,\g}z)-T^{(2)}_{\a\b,\g}(x*_{\a,\b}y,z)\\
&&+((123)(123)\o_H \hbox{id})T^{(2)}_{\g,\a}(z,x)*_{\g\a,\b}y-((123)(132)\o_H \hbox{id})(x*_{\a,\b\g}T^{(2)}_{\b,\g}(y,z))\\
&&+((123)\o_H \hbox{id})T^{(2)}_{\b,\g\a}(y, z*_{\g,\a}x)-((123)\o_H \hbox{id})T^{(2)}_{\b\g,\a}(y*_{\b,\g}z,x))\\
&&+((132)(123)\o_H \hbox{id})T^{(2)}_{\a,\b}(x,y)*_{\a\b,\g}z-((132)\o_H \hbox{id})T^{(2)}_{\g,\a}(z,x)*_{\g\a,\b}y\\
&&+((132)\o_H \hbox{id})T^{(2)}_{\g,\a\b}(z, x*_{\a,\b}y)\\
&&-((132)\o_H \hbox{id})T^{(2)}_{\g\a,\b}(z*_{\g,\a}x,y))\quad(\hbox{by Eqs. \eqref{L8}-\eqref{L9}})\\
&=&T^{(2)}_{\a,\b\g}(x, y*_{\b,\g}z)-T^{(2)}_{\a\b,\g}(x*_{\a,\b}y,z)+((123)\o_H \hbox{id})(T^{(2)}_{\b,\g\a}(y,z*_{\g,\a}x)\\
&&-T^{(2)}_{\b\g,\a}(y*_{\b,\g}z,x))+((132)\o_H \hbox{id})(T^{(2)}_{\g,\a\b}(z,x*_{\a,\b}y)-T^{(2)}_{\g\a,\b}(z*_{\g,\a}x,y))\\
&=&T^{(2)}\hat{\circ} T^{(0)}(x,y,z)_{\a,\b,\g}+((123)\o_H \hbox{id})T^{(2)}\hat{\circ} T^{(0)}(y,z,x)_{\b,\g,\a}\\
&&+((132)\o_H \hbox{id})T^{(2)}\hat{\circ} T^{(0)}(z,x,y)_{\g,\a,\b},
\end{eqnarray*}
finishing the proof.
$\hfill \square$
\\

{\bf Lemma 6.6} Under the same hypotheses as Lemma 6.5, the following identity holds for all $x, y, z\in A$ and $\a,\b,\g\in \Omega$:
\begin{eqnarray*}
&&((12)\o_H \hbox{id})T^{(1)}\hat{\circ} T^{(1)}(y,x,z)_{\b,\a,\g}+((13)\o_H \hbox{id})T^{(1)}\hat{\circ} T^{(1)}(z,y,x)_{\g,\b,\a}\\
&&+((23)\o_H \hbox{id})T^{(1)}\hat{\circ} T^{(1)}(x,z,y)_{\a,\g,\b}\\
&=&-((12)\o_H \hbox{id})T^{(2)}\hat{\circ} T^{(0)}(y,x,z)_{\b,\a,\g}-((13)\o_H \hbox{id})T^{(2)}\hat{\circ} T^{(0)}(z,y,x)_{\g,\b,\a}\\
&&-((23)\o_H \hbox{id})T^{(2)}\hat{\circ} T^{(0)}(x,z,y)_{\a,\g,\b}.
\end{eqnarray*}

{\bf Proof.} For all $x, y, z\in A$ and $\a,\b,\g\in \Omega$, we have
\begin{eqnarray*}
&&((12)\o_H \hbox{id})T^{(1)}\hat{\circ} T^{(1)}(y,x,z)_{\b,\a,\g}+((13)\o_H \hbox{id})T^{(1)}\hat{\circ} T^{(1)}(z,y,x)_{\g,\b,\a}\\
&&+((23)\o_H \hbox{id})T^{(1)}\hat{\circ} T^{(1)}(x,z,y)_{\a,\g,\b}\\
&=&-((12)\o_H \hbox{id})(T^{(0)}\hat{\circ} T^{(2)}(y, x, z)_{\b,\a,\g}+T^{(2)}\hat{\circ} T^{(0)}(y, x, z)_{\b,\a,\g})\\
&&-((13)\o_H \hbox{id})(T^{(0)}\hat{\circ} T^{(2)}(z, y, x)_{\g,\b,\a}+T^{(2)}\hat{\circ} T^{(0)}(z, y, x)_{\g,\b,\a})\\
&&-((23)\o_H \hbox{id})(T^{(0)}\hat{\circ} T^{(2)}(x, z, y)_{\a,\g,\b}+T^{(2)}\hat{\circ} T^{(0)}(x, z, y)_{\a,\g,\b})\\
&=&((12)\o_H \hbox{id})(T^{(2)}_{\b,\a}(y,x)*_{\b\a,\g}z-y*_{\b,\a\g}T^{(2)}_{\a,\g}(x,z)+T^{(2)}_{\b\a,\g}(y*_{\b,\a}x,z)\\
&&-T^{(2)}_{\b,\a\g}(y,x*_{\a,\g}z))+((13)\o_H \hbox{id})(T^{(2)}_{\g,\b}(z,y)*_{\g\b,\a}x-z*_{\g,\b\a}T^{(2)}_{\b,\a}(y,x)\\
&&+T^{(2)}_{\g\b,\a}(z*_{\g,\b}y,x)-T^{(2)}_{\g,\b\a}(z,y*_{\b,\a}x))+((23)\o_H \hbox{id})(T^{(2)}_{\a,\g}(x,z)*_{\a\g,\b}y\\
&&-x*_{\a,\g\b}T^{(2)}_{\g,\b}(z,y)+T^{(2)}_{\a\g,\b}(x*_{\a,\g}z,y)-T^{(2)}_{\a,\g\b}(x,z*_{\g,\b}y))\\
&=&((12)(132)\o_H \hbox{id})z*_{\g,\b\a}T^{(2)}_{\b,\a}(y,x)-((12)(123)\o_H \hbox{id})T^{(2)}_{\a,\g}(x,z)*_{\a\g,\b}y\\
&&+((12)\o_H \hbox{id})T^{(2)}_{\b\a,\g}(y*_{\b,\a}x,z)-((12)\o_H \hbox{id})T^{(2)}_{\b,\a\g}(y,x*_{\a,\g}z)\\
&&+((13)\o_H \hbox{id})T^{(2)}_{\g,\b}(z,y)*_{\g\b,\a}x-((13)\o_H \hbox{id})z*_{\g,\b\a}T^{(2)}(y,x)\\
&&+((13)\o_H \hbox{id})T^{(2)}_{\g\b,\a}(z*_{\g,\b}y,x)-((13)\o_H \hbox{id})T^{(2)}_{\g,\b\a}(z,y*_{\b,\a}x)\\
&&+((23)\o_H \hbox{id})T^{(2)}_{\a,\g}(x,z)*_{\a\g,\b}y-((23)(123)\o_H \hbox{id})T^{(2)}_{\g,\b}(z,y)*_{\g\b,\a}x\\
&&+((23)\o_H \hbox{id})T^{(2)}_{\a\g,\b}(x*_{\a,\g}z,y)-((23)\o_H \hbox{id})T^{(2)}_{\a,\g\b}(x,z*_{\g,\b}y)\quad\hbox{by Eqs. \eqref{L8}-\eqref{L9}}\\
&=&((12)\o_H \hbox{id})T^{(2)}_{\b\a,\g}(y*_{\b,\a}x,z)-((12)\o_H \hbox{id})T^{(2)}_{\b,\a\g}(y,x*_{\a,\g}z)\\
&&+((13)\o_H \hbox{id})T^{(2)}_{\g\b,\a}(z*_{\g,\b}y,x)-((13)\o_H \hbox{id})T^{(2)}_{\g,\b\a}(z,y*_{\b,\a}x)\\
&&+((23)\o_H \hbox{id})T^{(2)}_{\a\g,\b}(x*_{\a,\g}z,y)-((23)\o_H \hbox{id})T^{(2)}_{\a,\g\b}(x,z*_{\g,\b}y)\\
&=&-((12)\o_H \hbox{id})T^{(2)}\hat{\circ} T^{(0)}(y,x,z)_{\b,\a,\g}-((13)\o_H \hbox{id})T^{(2)}\hat{\circ} T^{(0)}(z,y,x)_{\g,\b,\a}\\
&&-((23)\o_H \hbox{id})T^{(2)}\hat{\circ} T^{(0)}(x,z,y)_{\a,\g,\b},
\end{eqnarray*}
finishing the proof.
$\hfill \square$
\\

In what follows, we give the main result of this section.
\\

{\bf Theorem 6.7} Let $\mathfrak{T}^{t}_{\a,\b}=\sum_{i=0}^{\infty}T^{(i)}_{\a,\b}t^{i}$ be the deformation of a commutative $\Omega$-associative $H$-pseudoalgebra $(A, (*_{\a,\b})_{\a,\b\in \Omega})$. Define $\{x*_{\a,\b}y\}=T^{(1)}_{\a,\b}(x,y)-(\sigma\o_H \hbox{id})T^{(1)}_{\b,\a}(y,x)$, for all $x, y\in A$ and $\a,\b\in \Omega$. Then $(A, *_{\a,\b}, \{\c*_{\a,\b}\c\})$ is an $\Omega$-Poisson $H$-pseudoalgebra.

{\bf Proof.} By Lemma 6.4, the compatible condition \eqref{6.1} holds. It is easy to prove that the pseudobracket $\{\cdot*_{\a,\b}\cdot\}$ is skew-commutative. It suffices to prove that $\{\cdot*_{\a,\b}\cdot\}$ satisfies the Jacobi identity. Using Lemmas 6.5 and 6.6, for all $x, y, z\in A$ and $\a,\b,\g\in \Omega$, we have
\begin{eqnarray*}
&&\{\{x*_{\a,\b}y\}*_{\a\b,\g}z\}-\{x*_{\a,\b\g}\{y*_{\b,\g}z\}\}+((12)\o_H \hbox{id})\{y*_{\b,\a\g}\{x*_{\a,\g}z\}\}\\
&=&T^{(1)}_{\a\b,\g}(T^{(1)}_{\a,\b}(x, y),z)-((132)\o_H \hbox{id})T^{(1)}_{\g,\a\b}(z,T^{(1)}_{\a,\b}(x,y))\\
&&-((12)\o_H \hbox{id})T^{(1)}_{\b\a,\g}(T^{(1)}_{\b,\a}(y,x),z)+((13)\o_H \hbox{id})T^{(1)}_{\g,\b\a}(z,T^{(1)}_{\b,\a}(y,x))\\
&&-T^{(1)}_{\a,\b\g}(x,T^{(1)}_{\b,\g}(y,z))+((123)\o_H \hbox{id})T^{(1)}_{\b\g,\a}(T^{(1)}_{\b,\g}(y,z),x)\\
&&+((23)\o_H \hbox{id})T^{(1)}_{\a,\g\b}(x,T^{(1)}_{\g,\b}(z,y))-((13)\o_H \hbox{id})T^{(1)}_{\g\b,\a}(T^{(1)}_{\g,\b}(z,y),x)\\
&&+((12)\o_H \hbox{id})T^{(1)}_{\b,\a\g}(y,T^{(1)}_{\a,\g}(x,z))-((23)\o_H \hbox{id})T^{(1)}_{\a\g,\b}(T^{(1)}_{\a,\g}(x,z),y)\\
&&-((123)\o_H \hbox{id})T^{(1)}_{\b,\g\a}(y,T^{(1)}_{\g,\a}(z,x))+((132)\o_H \hbox{id})T^{(1)}_{\g\a,\b}(T^{(1)}_{\g,\a}(z,x),y)\\
&=&T^{(1)}_{\a\b,\g}(T^{(1)}_{\a,\b}(x, y),z)-T^{(1)}_{\a,\b\g}(x,T^{(1)}_{\b,\g}(y,z))\\
&&-((123)\o_H \hbox{id})(T^{(1)}_{\b,\g\a}(y,T^{(1)}_{\g,\a}(z,x))-T^{(1)}_{\b\g,\a}(T^{(1)}_{\b,\g}(y, z),x))\\
&&-((132)\o_H \hbox{id})(T^{(1)}_{\g,\a\b}(z,T^{(1)}_{\a,\b}(x,y))-T^{(1)}_{\g\a,\b}(T^{(1)}_{\g,\a}(z, x),y))\\
&&+((12)\o_H \hbox{id})(T^{(1)}_{\b,\a\g}(y,T^{(1)}_{\a,\g}(x,z))-T^{(1)}_{\b\a,\g}(T^{(1)}_{\b,\a}(y, x),z))\\
&&+((13)\o_H \hbox{id})(T^{(1)}_{\g,\b\a}(z,T^{(1)}_{\b,\a}(y,x))-T^{(1)}_{\g\b,\a}(T^{(1)}_{\g,\b}(z, y),x))\\
&&+((23)\o_H \hbox{id})(T^{(1)}_{\a,\g\b}(x,T^{(1)}_{\g,\b}(z,y))-T^{(1)}_{\a\g,\b}(T^{(1)}_{\a,\g}(x, z),y))\\
&=&-T^{(1)}\hat{\circ} T^{(1)}(x,y,z)_{\a,\b,\g}-((123)\o_H \hbox{id})T^{(1)}\hat{\circ} T^{(1)}(y,z,x)_{\b,\g,\a}\\
&&-((132)\o_H \hbox{id})T^{(1)}\hat{\circ} T^{(1)}(z,x,y)_{\g,\a,\b}+((12)\o_H \hbox{id})T^{(1)}\hat{\circ} T^{(1)}(y,x,z)_{\b,\a,\g}\\
&&+((13)\o_H \hbox{id})T^{(1)}\hat{\circ} T^{(1)}(z,y,x)_{\g,\b,\a}+((23)\o_H \hbox{id})T^{(1)}\hat{\circ} T^{(1)}(x,z,y)_{\a,\g,\b}\\
&=&T^{(2)}\hat{\circ} T^{(0)}(x,y,z)_{\a,\b,\g}+((123)\o_H \hbox{id})T^{(2)}\hat{\circ} T^{(0)}(y,z,x)_{\b,\g,\a}\\
&&+((132)\o_H \hbox{id})T^{(2)}\hat{\circ} T^{(0)}(z,x,y)_{\g,\a,\b}-((12)\o_H \hbox{id})T^{(2)}\hat{\circ} T^{(0)}(y,x,z)_{\b,\a,\g}\\
&&-((13)\o_H \hbox{id})T^{(2)}\hat{\circ} T^{(0)}(z,y,x)_{\g,\b,\a}-((23)\o_H \hbox{id})T^{(2)}\hat{\circ} T^{(0)}(x,z,y)_{\a,\g,\b}\\
&=&T^{(2)}_{\a,\b\g}(x,y*_{\b,\g}z)-T^{(2)}_{\a\b,\g}(x*_{\a,\b}y,z)\\
&&+((123)\o_H \hbox{id})T^{(2)}_{\b,\g\a}(y,z*_{\g,\a}x)-((123)\o_H \hbox{id})T^{(2)}_{\b\g,\a}(y*_{\b,\g}z,x)\\
&&+((132)\o_H \hbox{id})T^{(2)}_{\g,\a\b}(z,x*_{\a,\b}y)-((132)\o_H \hbox{id})T^{(2)}_{\g\a,\b}(z*_{\g,\a}x,y)\\
&&-((12)\o_H \hbox{id})T^{(2)}_{\b,\a\g}(y,x*_{\a,\g}z)+((12)\o_H \hbox{id})T^{(2)}_{\b\a,\g}(y*_{\b,\a}x,z)\\
&&-((13)\o_H \hbox{id})T^{(2)}_{\g,\b\a}(z,y*_{\b,\a}x)+((13)\o_H \hbox{id})T^{(2)}_{\g\b,\a}(z*_{\g,\b}y,x)\\
&&-((23)\o_H \hbox{id})T^{(2)}_{\a,\g\b}(x,z*_{\g,\b}y)+((23)\o_H \hbox{id})T^{(2)}_{\a\g,\b}(x*_{\a,\g}z,y)\\
&=&T^{(2)}_{\a,\b\g}(x,y*_{\b,\g}z)-T^{(2)}_{\a\b,\g}(x*_{\a,\b}y,z)\\
&&+((123)(23)\o_H \hbox{id})T^{(2)}_{\b,\a\g}(y,x*_{\a,\g}z)-((123)\o_H \hbox{id})T^{(2)}_{\b\g,\a}(y*_{\b,\g}z,x)\\
&&+((132)\o_H \hbox{id})T^{(2)}_{\g,\a\b}(z,x*_{\a,\b}y)-((132)(12)\o_H \hbox{id})T^{(2)}_{\a\g,\b}(x*_{\a,\g}z,y)\\
&&-((12)\o_H \hbox{id})T^{(2)}_{\b,\a\g}(y,x*_{\a,\g}z)+((12)\o_H \hbox{id})T^{(2)}_{\b\a,\g}(y*_{\b,\a}x,z)\\
&&-((13)(23)\o_H \hbox{id})T^{(2)}_{\g,\a\b}(z,x*_{\a,\b}y)+((13)(12)\o_H \hbox{id})T^{(2)}_{\b\g,\a}(y*_{\b,\g}z,x)\\
&&-((23)\o_H \hbox{id})T^{(2)}_{\a,\g\b}(x,z*_{\g,\b}y)+((23)\o_H \hbox{id})T^{(2)}_{\a\g,\b}(x*_{\a,\g}z,y)\quad\hbox{by Eqs. \eqref{L10}-\eqref{L11}}\\
&=&0.
\end{eqnarray*}
This completes the proof.
$\hfill \square$
\\

\section*{Acknowledgements}

This work was supported by the National Natural Science Foundation of China (Nos. 12301023, 12201188).

\section*{Date Availability}

Date sharing is not applicable to this article as no new date were created or analyzed in
this study.

\end{document}